\documentclass[a4paper,12pt]{article}
\usepackage[centertags]{amsmath}
\usepackage[english]{babel}
\usepackage{amsfonts}
\usepackage{amsmath,amssymb,graphics,amscd}
\usepackage{amsthm}
\usepackage[all]{xy}
\newcommand{\semi}{{\,\rule[.1pt]{.4pt}{5.3pt}\hskip-1.9pt\times}}

\newcommand{\al}{\alpha}

\newcommand{\avee}{\alpha^\vee}
\newcommand{\bvee}{\beta^\vee}

\newcommand{\eps}{\epsilon}
\newcommand{\om}{\omega}

\newcommand{\lam}{{\lambda}}
\newcommand{\Lam}{{\Lambda}}

\newcommand{\Lhg}{\widehat{\mathfrak{g}}}
\newcommand{\Lhh}{\widehat{\mathfrak{h}}}
\newcommand{\Lhb}{\widehat{\mathfrak{b}}}
\newcommand{\Lhn}{\widehat{\mathfrak{n}}}
\newcommand{\wDelta}{\widehat{{\Delta}}}
\newcommand{\wPhi}{\widehat{{\Phi}}}
\newcommand{\ph}{\widehat{P}}

\newcommand{\waff}{{W^{\rm aff}}}
\newcommand{\wiff}{{\widetilde{W}^{\rm aff}}}
\newcommand{\ca}{\mathcal{A}}

\newcommand{\chp}{{\check{P}}}

\newcommand{\semc}{{\rm sc}}
\newcommand{\Hom}{{{\hbox{\rm Hom}}\,}}

\newcommand{\bn}{{\mathbb N}}
\newcommand{\br}{{\mathbb R}}
\newcommand{\bz}{{\mathbb Z}}

\newcommand{\bc}{{\mathbb C}}

\newcommand{\Lg}{{\mathfrak{g}}}
\newcommand{\Lb}{{\mathfrak{b}}}
\newcommand{\Lh}{{\mathfrak{h}}}
\newcommand{\Ln}{{\mathfrak{n}}}

\newcommand{\Lgl}{{\mathcal{L}\mathfrak{g}}}
\newcommand{\Lgc}{{\mathcal{C}\mathfrak{g}}}
\newcommand{\Lbc}{{\mathfrak{I_b}}}
\newtheorem{thm}{Theorem}
\newtheorem*{thm*}{Theorem}
\newtheorem{defn}{Definition}

\newtheorem{conj}{Conjecture}
\newtheorem{rem}{Remark}

\newtheorem{prop}{Proposition}
\newtheorem{lem}{Lemma}
\newtheorem{coro}{Corollary}
\def\proof{\noindent{\it Proof. }}
\def\endpf{\hfill$\bullet$\vskip 5pt}
\font\eightsc=cmcsc10 scaled800

\begin{document}
\title{Weyl modules, Demazure modules, KR-modules,
crystals, fusion products and limit constructions}
\author{G. Fourier$^*$ and P. Littelmann\footnote{This research has been partially
supported by the EC TMR network "LieGrits", contract
MRTN-CT 2003-505078 and the DFG-Graduiertenkolleg 1052,
{\sl 2000 Mathematics Subject Classification. 22E46, 14M15.}}}
\maketitle
\begin{abstract}
We study finite dimensional representations
of current algebras, loop algebras and their quantized versions. 
For the current algebra of a simple Lie algebra of type {\tt ADE}, 
we show that Kirillov-Reshetikhin modules and Weyl modules are in fact all 
Demazure modules. As a consequence one obtains an elementary 
proof of the dimension formula for Weyl modules for the current and the 
loop algebra. Further, we show that the crystals of the Weyl 
and the Demazure module are the same up to some additional label zero 
arrows for the Weyl module. 

For the current algebra $\Lgc$ of an arbitrary simple Lie algebra, 
the fusion product of Demazure modules of the same level
turns out to be again a Demazure module. As an application we 
construct the $\Lgc$-module structure of the Kac-Moody algebra 
$\Lhg$-module $V(\ell\Lam_0)$ as a semi-infinite fusion product of
finite dimensional $\Lgc$--modules. 
\end{abstract}
\tableofcontents%
%%%%%%%%%%%%%%%%%%%%%%%%%%%%%%%%%%%%%%%%%%%%%%%%%%%%%%%%%%%%%%%%%%%%%%%%%%%%%%%%%%%%%%%%%%%%%%%%%%%%%%%%%%%%%%%%%%%%%%%%%%%%%%%%%%%
%         Introduction
%%%%%%%%%%%%%%%%%%%%%%%%%%%%%%%%%%%%%%%%%%%%%%%%%%%%%%%%%%%%%%%%%%%%%%%%%%%%%%%%%%%%%%%%%%%%%%%%%%%%%%%%%%%%%%%%%%%%%%%%%%%%%%%%%%%
\section{Introduction}
Let $\Lg$ be a semisimple complex Lie algebra.
The theory of finite dimensional representations of
its loop algebra $\Lgl=\Lg\otimes\bc[t,t^{-1}]$, its quantized loop algebra $U_q(\Lgl)$ and
its current algebra  $\Lgc=\Lg\otimes\bc[t]$ have been the subject of many articles
in the recent years.
See for example \cite{AK}, \cite{BN}, \cite{C1}, \cite{CL}, \cite{CM1}, \cite{CM3},
\cite{CP1}, \cite{CP2}, \cite{CP3}, \cite{FKL}, \cite{FM}, \cite{FR},
\cite{Ka2}, \cite{Ke}, \cite{N} for different approaches and
different aspects of this subject.

The notion of a Weyl module in this context was introduced in \cite{CP1}
for the affine Kac-Moody algebra and its quantized version. These
modules can be described in terms of generators and relations,
and they are characterized by the following universal property:
any finite dimensional highest weight module which is generated
by a one dimensional highest weight space, is a quotient of
a Weyl module. This notion can be naturally extended
to the category of finite dimensional representations of the
current algebra (\cite{CL}, \cite{FL2}).
Another intensively studied class of modules are the Kirillov-Reshetikhin
modules, a name that originally refers to evaluation modules of the Yangian.
In \cite{C2} Chari gave a definition of these modules for the current algebra
in terms of generators and relations.

The current algebra is a subalgebra of a maximal parabolic subalgebra
of the affine Kac-Moody algebra  $\Lhg$. Let $\Lam$ be a dominant weight
for $\Lhg$ and denote by $V(\Lam)$ the associated
(infinite dimensional) irreducible $\Lhg$-representation.
Another natural class of finite dimensional representations of the
current algebra are provided by certain
Demazure submodules of $V(\Lam)$.
Of particular interest for this paper are the twisted
(see section~\ref{defndemazure}) $\Lgc$-stable Demazure
submodules $D(m,\lam)$ of $V(m\Lam_0)$, where $\Lam_0$ is the fundamental
weight associated  to the additional node of the extended Dynkin
diagram of $\Lg$.

If $\Lg$ is simply laced, then we can identify the weight and the coweight lattice,
so the Weyl modules as well as the twisted
$\Lgc$-stable Demazure submodules of $V(m\Lam_0)$ are classified
by dominant weights $\lam\in P^+$. 
\vskip 5pt\noindent
{\bf Theorem A}
{\it For a simple complex Lie algebra of simply laced type,
the Weyl module $W(\lam)$ and the Demazure module $D(1,\lam)$ are isomorphic as $\Lgc$-modules.}
\vskip 5pt
Also the Demazure modules of higher level are related to an interesting class
of finite dimensional modules for $\Lgc$. Let $\Lg$ be
an arbitrary simple complex Lie algebra, the  $\Lgc$-stable Demazure modules 
$D(m,\lam^\vee)$ are classified by dominant coweights $\lam\in \check{P}^+$. 
\vskip 5pt\noindent
{\bf Theorem B}
{\it For a fundamental coweight $\om_i^{\vee}$ let $d_i=1,2$ or $3$ be such that 
$d_i \om_i = \nu(\om_i^{\vee})$. The Kirillov-Reshetikhin module $KR(d_im\om_i)$
is, as $\Lgc$-module, isomorphic to the Demazure module $D(m, \om_i^{\vee})$.
In particular, in the simply laced case all Kirillov-Reshetikhin modules are Demazure 
modules.}
\begin{rem}\rm
The fact that $D(m, \om_i^{\vee})$ is a quotient of a  Kirillov-Reshetikhin module
has been already pointed out in \cite{CM2}. In the same paper Chari and Moura
have shown that $D(m, \om_i^{\vee})$ is isomorphic to $KR(d_im\om_i)$ for
all classical groups using character calculations. Our proof
is independent of the type of the algebra.
\end{rem}
To stay inside the class of cyclic highest weight modules, the 
tensor product of cyclic $\Lgc$-modules is often replaced by the fusion product
of modules \cite{FL1}. 
%We show that a Demazure module is actually the fusion product of ``smaller" Demazure modules
\vskip 5pt\noindent
{\bf Theorem C}
{\it Let $\Lg$ be a complex simple Lie algebra and let $\lam^\vee=\lam^\vee_1+\ldots+\lam^\vee_r$
be a decomposition of a dominant coweight as a sum of dominant coweights. Then
$D(m,\lam^\vee)$  and the fusion product $D(m,\lam_1^\vee)*\cdots*D(m,\lam_r^\vee)$ are isomorphic
as  $\Lgc$-modules.}
%As modules for the current algebra $\Lgc$, the Demazure module
%$D(m,\lam^\vee)$  and the fusion product 
%$D(m,\lam_1^\vee)*\cdots*D(m,\lam_r^\vee)$ are isomorphic.}
\begin{rem}\rm 
The theorem shows in particular that the fusion
product of Demazure modules of the same level is associative and independent 
of the parameters used in the fusion construction. In \cite{AK}
it is shown that the fusion product of Kirillov-Reshetikhin modules
of arbitrary levels is independent of the parameters. 
\end{rem}
As a consequence we obtain for the Weyl module $W(\lam)$ in the simply laced case:
\vskip 5pt\noindent
{\bf Corollary A}
{\it Suppose $\Lg$ is of simply laced type. Let $\lam=a_1\om_1+\ldots+a_n\om_n$ be a decomposition
of a dominant weight $\lam\in P^+$ as a sum of fundamental weights.
Then the Weyl module $W(\lam)$ for the current algebra
is the fusion product of the fundamental Weyl modules:
$$
W(\lam)\simeq
\underbrace{W(\om_1)*\cdots * W(\om_1)}_{a_1}
*\cdots *
\underbrace{W(\om_n)*\cdots * W(\om_n)}_{a_n}.
$$}
The Weyl modules for the loop algebra
are classified by $n$-tuples $\pi=(\pi_1,\ldots,\pi_n)$
of polynomials $\pi_j\in\bc[u]$ with constant term 1 \cite{CP1}.
The associated dominant weight is $\lam_\pi= \sum_i\deg\pi_i\om_i$.
Similarly, the Weyl modules for the quantized loop algebra
are classified by $n$-tuples $\pi_q=(\pi_{q,1},\ldots,\pi_{q,n})$
of polynomials $\pi_{q,j}\in\bc(q)[u]$ with constant term 1, the associated
weight $\lam_{\pi_q}$ is defined as above. 

It was conjectured in \cite{CP1} (and proved in the ${\mathfrak{sl}}_2$-case)
that the dimension of the Weyl modules depend only  $\lam_{\pi_q}$ respectively $\lam_\pi$.
More precisely, they conjectured that the dimension is the (appropriate) 
product of the dimension of the ``fundamental modules''.
As Hiraku Nakajima has pointed out to us, the dimension conjecture can be deduced 
using the theory of global basis. 
The results of Kashiwara \cite{Ka2,Ka3} %12, 13] 
imply that the Weyl modules are specializations (the q = 1 limit) of certain 
finite-dimensional quotients of the extremal weight modules for the quantum affine algebra. 
The results of Beck and Nakajima \cite{BN, N,N2} imply 
that these quotients (and hence their specializations) have the correct dimension.

A different approach to prove the dimension conjecture was 
suggested in \cite{CP1,CP2}. In fact, using the specialization and 
dimension arguments outlined there, in the simply laced case 
the dimension formula is an immediate
consequence of Theorem A and Theorem C:
\vskip 5pt\noindent
{\bf Corollary B}
{\it Let $\Lg$ be a simple Lie algebra of simply laced type, 
let $\lam=\sum m_i\om_i$ be a dominant
weight (for $\Lg$), let $\pi$ (resp. $\pi_q$) be an $n$-tuple of polynomials in $\bc[u]$ (resp. in $\bc(q)[u]$) 
with constant term 1 such that $\lam=\lam_\pi=\lam_{\pi_q}$. Then}
$$
\dim W(\lam) = \dim W(\pi) = \dim W_{q}(\pi_{q}) = \dim D(1,\lam^{\vee}) = \prod_i (\dim W(\om_i))^{m_i}
$$
\begin{rem}\rm
For $\Lg={\mathfrak{sl}_n}$, the
connection between Demazure modules in $V(\Lam_0)$ and 
Weyl modules had been already obtained by Chari and Loktev in \cite{CL}. 
The isomorphism between the Weyl module $W(\lam)$ and the 
Demazure module $D(1,\lam)$ has been conjectured in \cite{FKL}.
\end{rem}
For a dominant weight $\lam=\sum m_i\om_i$ let $\pi_{\lam,a}$ be the
tuple having $(1-au)^{m_i}$ as $i$-th entry. The quantum Demazure
module $D_q(m,\lam)$ has an associated crystal graph which is a subgraph
of the crystal graph of the corresponding irreducible $U_q(\Lhg)$-representation.
We conjecture that by adding appropriate label zero arrows,
one gets the graph of an irreducible $U_q(\Lhg)$-representation.
In the simply laced case and level one we have:
\vskip 5pt\noindent
{\bf Proposition\/}
{\it The crystal graph of $D_q(1,\lam)$ is obtained from the crystal
graph of $W_q(\pi_{q,\lam,1})$ by omitting certain label zero arrows.
More precisely, let $B(\lam)_{cl}$ be the path model for
$W_q(\pi_{q,\lam,1})$ described in \cite{NS3}, then the crystal
graph of the Demazure module is isomorphic to the graph of the concatenation 
$\pi_{\Lambda_0}*B(\lam)_{cl}$.}
%$D_q(1,\lam)$ admits the structure of an 
%irreducible $U_q(\Lgl)$-module, isomorphic to $W_q(\pi_{q,\lam,1})$, with crystal basis, 
%such that the $U_q(\Lg)$-module structure coincides with the natural $U_q(\Lhg)$-structure
%coming from the Demazure module construction. The crystal
%graph of the Demazure module is isomorphic to the graph of the concatenation 
%$\pi_{\Lambda_0}*B(\lam)_{cl}$, where $B(\lam)_{cl}$ is the path model for
%$W_q(\pi_{q,\lam,1})$ described in \cite{NS3}.
\vskip 5pt
%We conjecture (see \cite{FoL}) that the proposition above holds 
%in general for arbitrary level,
%but the Weyl modules should be replaced by irreducible modules. 
In the simply laced case, the restriction of the loop Weyl module $W(\pi_{\lam,a})$ 
to $\Lgc$ is (up to a twist by an automorphism) the Weyl module $W(\lam)$. It follows:
\vskip 5pt\noindent
{\bf Corollary C}
{\it The Demazure module $D(m, \lam)$ of level $m$ can be equipped with the structure of 
a cyclic $U(\Lgl)$-module such that the $\Lg$-module structure coincides with the natural 
$\Lg$-structure coming from the Demazure module construction.}
\vskip 5pt
Let $V(m\Lam_0)$, $m\in\bn$, be the irreducible
highest weight module of highest weight $m\Lam_0$
for the affine Kac-Moody $\Lhg$. In \cite{FoL}
we gave a description of the $\Lg$-module structure
of this representation in terms of a semi-infinite tensor product.
Using Theorem~C, we are able to lift this result to the level of modules
for the current algebra. The theorem holds in a much more general setting (see
Remark~\ref{moregenerallimit}), but for the convenience of a uniform
presentation, let $\Theta$ be
the highest root of the root system of $\Lg$.
\vskip 5pt\noindent
{\bf Theorem D}
{\it Let $D(m, n\Theta)\subset V(m\Lam_0)$ be the
Demazure module of level $m$ corresponding to the
translation at $-n\Theta$. Let $w\neq 0$ be a $\Lgc$-invariant
vector of $D(m,\Theta)$. Let $V^{\infty}_{m} $ be the direct limit
$$
D(m, \Theta) \hookrightarrow D(m, \Theta) \ast D(m,\Theta) \hookrightarrow D(m,\Theta) \ast 
D(m,\Theta) \ast D(m,\Theta) \hookrightarrow  \ldots
$$
where the inclusions are given by $v \mapsto w \otimes v$. \\
Then $V(m\Lam_0)$ and $V^{\infty}_{m}$ are isomorphic as $U(\Lgc)$-modules.}
\vskip 5pt
The semi-infinite fusion construction can be seen as an extension of the 
construction of Feigin and Feigin \cite{FF} ($\Lg={\mathfrak{sl}_2}$) and
Kedem \cite{Ke} ($\Lg={\mathfrak{sl}_n}$) to arbitrary simple Lie algebras.
We conjecture (see Conjecture~\ref{semiinfconjecture}) 
that, as in \cite{FF} and \cite{Ke}, the semi-infinite
fusion construction works for arbitrary dominant weights and not
only for multiples of $\Lam_0$. %\marginpar{ge\"andert}
\begin{rem}\rm
Naito and Sagaki \cite{NS1}, \cite{NS2}, \cite{NS3} gave a path model for the 
Weyl modules $W(\om)$ for all fundamental weights and $\Lg$ of arbitrary type.  
Since the Weyl modules coincide with the level-one Demazure modules 
provided $\Lg$ is simply-laced, the semi-infinite limit construction above
gives on the combinatorial side a combinatorial limit path model for the 
representation $V(\Lam_0)$ as a semi-infinite concatenation of a finite 
path model, extending in this sense the approach of Magyar in \cite{Mag}.
\end{rem}

After introducing some notation, we will recall in more
detail the definition of Demazure and Weyl modules and fusions products.
The proof of the Theorems~A -- C, their corollaries and the proposition
is given in section~\ref{connections} (see Theorem~ \ref{KRgleichD}, 
\ref{WeylgleichDemazure} and \ref{demazure-zerlegung}). 
The proof of Theorem~D is given in section~\ref{seclimitconstruction},
see Theorem~\ref{limitthm}.

{\bf Acknowledgements.} We are grateful to V. Chari, A. Joseph,
S.~Loktev, T.~Miwa, and A.~Moura for many helpful discussions and
useful hints. We would like to thank H. Nakajima for the discussions 
concerning the dimension conjecture.
%%%%%%%%%%%%%%%%%%%%%%%%%%%%%%%%%%%%%%%%%%%%%%%%%%
%%%%%%%%%%%%%%%%%%%%%%%%%%%%%%%%%%%%%%%%%%%%%%%%%%
%%%%%%%%%%%%%%%%%%%%%%%%%%%%%%%
%          Notation and basics
%%%%%%%%%%%%%%%%%%%%%%%%%%%%%%%%%%%%%%%%%%%%%%%%%%
%%%%%%%%%%%%%%%%%%%%%%%%%%%%%%%%%%%%%%%%%%%%%%%%%%
%%%%%%%%%%%%%%%%%%%%%%%%%%%%%%%
\section{Notation and basics}
\subsection{Affine Kac-Moody algebras}
In this section we fix the notation and the usual technical padding.
Let $\Lg$ be simple complex Lie algebra. We fix a Cartan subalgebra
$\Lh$ in $\Lg$ and a Borel subalgebra $\Lb\supseteq \Lh$. Denote
$\Phi\subseteq \Lh^*$ the root system of $\Lg$, and, corresponding to the choice
of $\Lb$, let $\Phi^+$ be the set of positive roots
and let $\Delta=\{\al_1,\ldots,\al_n\}$ be the corresponding basis of $\Phi$.

For a root $\beta\in\Phi$ let $\beta^\vee\in\Lh$ be its coroot.
The basis of the dual root system (also called the coroot system)
$\Phi^\vee\subset \Lh$ is denoted $\Delta^\vee=\{\al_1^\vee,\ldots,\al_n^\vee\}$.

We denote throughout the paper by $\Theta=\sum_{i=1}^n a_i\al_i$ the
highest root of $\Phi$ and by $\Theta^\vee=\sum_{i=1}^n a_i^\vee\al_i^\vee$
its coroot. Note that $\Theta^\vee$ is in general not the highest root of
$\Phi^\vee$. (For more details concerning the connection
with the dual root system of the affine root system $\wPhi$ see \cite{Kac},
Chapter 6.)

The Weyl group $W$ of $\Phi$ is generated
by the simple reflections $s_i=s_{\al_i}$ associated to the simple roots.

Let $P$ be the weight lattice of $\Lg$ and let $P^+$ be the subset
of dominant weights.  The group algebra of $P$ is denoted $\bz[P]$,
we write $\chi=\sum  m_\mu e^\mu$ (finite sum, $\mu\in P$, $m_\mu\in\bz$)
for an element in $\bz[P]$, where the embedding $P\hookrightarrow \bz[P]$ is
defined by $\mu\mapsto e^\mu$.

We denote the coweight lattice by $\chp$, i.e.,
this is the lattice of integral weights for the dual root root system.
The dominant coweights are denoted $\chp^+$.

Corresponding to the enumeration of the simple roots
let $\om_1,\ldots,\om_n$ be the fundamental weights.
Let $\Lh_\br$ be the ``real part'' of $\Lh$, i.e., $\Lh_\br$ is
the real span in $\Lh$ of the coroots $\al_1^\vee,\ldots,\al_n^\vee$,
and let $\Lh_\br^*$ be the real span of the fundamental weights
$\om_1,\ldots,\om_n$. Let $(\cdot,\cdot)$ be the unique invariant
symmetric non-degenerate bilinear form on $\Lg$ normalized
such that the restriction to $\Lh$ induces an  isomorphism
$$
\nu:\Lh_\br\longrightarrow \Lh_\br^*, \quad
\nu(h):\left\{
\begin{array}{rcl}
\Lh&\rightarrow&\br\\
h'&\mapsto &(h,h')\\
\end{array}\right.
$$
mapping $\Theta^\vee$ to $\Theta$. 
With the notation as above it follows
for the weight lattice $P^\vee$ of the dual root system $\Phi^\vee$ that
$$
\nu(\al_i^\vee)=\frac{a_i}{a_i^\vee}\al_i\quad{\rm and}\quad
\nu(\om_i^\vee)=\frac{a_i}{a_i^\vee}\om_i,\quad\forall\,i=1,\ldots,n.
$$
Let $\Lhg$ be the affine Kac--Moody algebra corresponding
to the extended Dynkin diagram of $\Lg$ (see \cite{Kac}, Chapter 7):
$$
\Lhg=\Lg\otimes_\bc\bc[t,t^{-1}]\oplus \bc K\oplus \bc d
$$
Here $d$ denotes the derivation $d=t\frac{d}{dt}$, $K$ is
the canonical central element, and the Lie bracket is given by
\begin{equation}\label{Lieklammer}
[t^m\otimes x + \lam K +\mu d, t^n\otimes y + \nu K + \eta d]
=t^{m+n}\otimes [x,y] + \mu n t^n\otimes y +\eta m t^m\otimes x +
m\delta_{m,-n}(x, y) K.
\end{equation}
The Lie algebra $\Lg$ is
naturally a subalgebra of  $\Lhg$.
In the same way, the Cartan subalgebra $\Lh\subset \Lg$ and the Borel
subalgebra $\Lb\subset \Lg$ are
subalgebras of the Cartan subalgebra $\Lhh$ respectively the Borel
subalgebra $\Lhb$ of $\Lhg$:
%%%%%%%%%%%%%%%%%%%%%%%%%%%%%%%%%%%%%%%%
\begin{equation}\label{hdecomp}
%%%%%%%%%%%%%%%%%%%%%%%%%%%%%%%%%%%%%%%%
\Lhh=\Lh\oplus\bc K\oplus\bc d,\quad
\Lhb=\Lb\oplus\bc K\oplus\bc d\oplus\Lg\otimes_\bc t\bc[t]
\end{equation}
Denote by $\widehat\Phi$ the root system of $\Lhg$ and
let $\wPhi^+$ be the subset of positive roots. The positive non-divisible
imaginary root in $\wPhi^+$ is denoted $\delta$. The simple roots
are $\wDelta=\{\al_0,\al_1,\ldots,\al_n\}$ where $\al_0=\delta-\Theta$.
We identify the root system $\Phi$ of $\Lg$ with the root subsystem
of $\wPhi$ generated by the simple roots $\al_1,\ldots,\al_n$.

Let $\Lam_0,\ldots,\Lam_n$ be the corresponding fundamental weights,
then for $i=1,\ldots,n$ we have
\begin{equation}\label{Lamdecomp}
\Lam_i=\om_i+a_i^\vee\Lam_0.
\end{equation}
The decomposition of $\Lhh$ in (\ref{hdecomp}) has its corresponding
version for the dual space $\Lhh^*$:
\begin{equation}\label{hveedecomp}
\Lhh^* = \Lh^*\oplus\bc\Lam_0\oplus\bc \delta
\end{equation}
Here the elements of
$\Lh^*$ are extended trivially, $\langle\Lambda_0,\Lh\rangle
=\langle\Lambda_0,d\rangle=0$ and $\langle\Lambda_0,K\rangle=1$,
and $\langle\delta,\Lh\rangle=\langle\delta,K\rangle=0$ and
$\langle\delta,d\rangle=1$. Let $\wDelta^\vee=\{\al_0^\vee,\al^\vee_1,
\ldots,\al^\vee_n\}\subset \Lhh$ be the corresponding basis of the
coroot system, then $\al_0^\vee=K-\Theta^\vee$.
Recall that the positive affine roots are precisely the roots of the form
$$
\wPhi^+=\{\beta+s\delta\mid\beta\in\Phi^+,\,s\ge 0\}\cup
\{-\beta+s\delta\mid\beta\in\Phi^+,\,s> 0\}\cup \{s\delta\mid s>0\}
$$
For a real positive root $\beta+s\delta$ respectively $-\beta+s\delta$
the corresponding coroot is
\begin{equation}\label{coroots}
(\beta+s\delta)^\vee=\beta^\vee +s\frac{(\beta^\vee,\beta^\vee)}{2}K
\quad\mbox{respectively}
\quad
(-\beta+s\delta)^\vee=-\beta^\vee+s\frac{(\beta^\vee,\beta^\vee)}{2}K
\end{equation}
Let $\Lhh_\br^*$ be the real span $\br\delta+\sum_{i=0}^n\br\Lam_i$,
note that by the decomposition (\ref{hveedecomp}) and by (\ref{Lamdecomp})
we have  $\Lh_\br^* \subseteq \Lhh_\br^*$.
The affine Weyl group  $\waff$ is generated by the
reflections $s_0, s_1, . . . , s_n$ acting on  $\Lhh_\br^*$. (We use again
the abbreviation $s_i=s_{\al_i}$ for a simple root $\al_i$.) The cone
$\widehat{C} = \{ \Lam\in\Lhh_\br^* |\langle\Lam,\al_i^\vee\rangle\ge 0,
i= 0, . . . , n\}$ is the  fundamental Weyl chamber for $\Lhg$.

We keep the convention and put a $\widehat\ $ on (almost) everything
related to $\Lhg$. We denote by $\ph$ the weight lattice of $\Lhg$ and by
$\ph^+$ the subset of dominant weights. As before, let
$\bz[\ph]$ be the group algebra of $\ph$, so an element
in the algebra is a finite sum of the form
$\sum  m_\mu e^\mu$,
$\mu\in \ph$ and $ m_\mu \in \bz$. Recall the following special properties
of the imaginary root $\delta$ (see for example \cite{Kac}, Chapter 6):
\begin{equation}\label{Kaczitat}
\langle\delta,\al_i^\vee\rangle=0\,\forall\, i=0,\ldots,n\quad   w(\delta)=\delta
\,\forall\, w\in \waff,\quad
\langle\al_0,\al_i^\vee\rangle=-\langle\Theta,\al_i^\vee\rangle\,{\rm for\ }i\ge 1
\end{equation}
Put $a_0=a_0^\vee=1$ and let $A=(a_{i,j})_{0\le i,j\le n}$ be the (generalized)
Cartan matrix of $\Lhg$. We have a non--degenerate symmetric bilinear
form $(\cdot,\cdot)$ on $\Lhh$ defined by (\cite{Kac}, Chapter 6)
\begin{equation}\label{hatform}
\left\{
\begin{array}{ll}
(\al_i^\vee,\al_j^\vee)=\frac{a_j}{a_j^\vee}a_{i,j}&i,j=0,\ldots,\ell\\
(\al_i^\vee,d)=0&i=1,\ldots,\ell\\
(\al_0^\vee,d)=1&(d,d)=0.\\
\end{array}\right.
\end{equation}
The corresponding isomorphism $\nu:\Lhh\rightarrow\Lhh^*$ maps
$$
\nu(\al_i^\vee)=\frac{a_i}{a_i^\vee}\al_i,\quad
\nu(K)=\delta,\quad \nu(d)=\Lambda_0.
$$\\
Since $\waff$ fixes $\delta$, the affine Weyl $\waff$ can be defined as the
subgroup of $GL(\Lh_{\semc,\br}^*)$
generated by the induced reflections $s_0,\ldots,s_n$. Another well--known
description of the affine Weyl group is the following. Let $M\subset \Lh^*_\br$
be the lattice $M=\nu(\bigoplus_{i=1}^n\bz\al_i^\vee)$. If $\Lg$ is simply
laced, then $M$ is the root lattice in $\Lh^*_\br$, otherwise $M$ is the lattice
in $\Lh^*_\br$ generated by the long roots.

An element $\Lam\in \Lh_{\semc,\br}^*$ can be uniquely decomposed
into $\Lam=\lam +b\Lam_0$ such that $\lam\in \Lh^*_\br$.
For an element $\mu\in M$ let $t_\mu\in GL(\Lh_{\semc,\br}^*)$ be the map
defined by
%%%%%%%%%%%%%%%%%%%%%%%%%%%%%%%%%%%%%%%%%%
\begin{equation}\label{muaction}
\Lam =\lam +b\Lam_0\mapsto
t_\mu(\Lam)=\lam +b\Lam_0 + b\mu=\Lam+\langle \Lam, K\rangle\mu.
\end{equation}
%%%%%%%%%%%%%%%%%%%%%%%%%%%%%%%%%%%%%%%%%%
Obviously we have $t_{\mu}\circ t_{\mu'}=t_{\mu+\mu'}$, denote $t_M$ the
abelian subgroup of $GL(\Lh_{\semc,\br}^*)$consisting of the elements
$t_\mu$, $\mu\in M$. Then $\waff$ is the semi-direct product
$\waff= W\semi t_M$.

The {\it extended affine Weyl group} $\wiff$ is the semi-direct product
$\wiff= W\semi t_L$, where $L=\nu(\bigoplus_{i=1}^n\bz\om^\vee_i)$
is the image of the coweight lattice. The action of an element $t_\mu$,
$\mu\in L$, is defined as above in (\ref{muaction}).

Let $\Sigma$ be the subgroup of $\wiff$ stabilizing the dominant Weyl chamber
$\widehat{C}$:
$$
\Sigma=\{\sigma\in\wiff\mid \sigma(\widehat{C})=\widehat{C}\}.
$$
Then $\Sigma$ provides a complete system of coset representatives of
$\wiff/\waff$, so we can write in fact $\wiff=\Sigma\semi\waff$.

The elements $\sigma\in\Sigma$ are all of the form
$$
\sigma=\tau_i t_{-\nu(\omega_i^\vee)}= \tau _i t_{-\omega_i},
$$
where $\omega_i^\vee$ is a minuscule fundamental coweight. Further,
set $\tau_i=w_0 w_{0,i}$, where $w_0$ is the longest word $W$ and
$w_{0,i}$ is the longest word in $W_{\om_i}$, the stabilizer of $\om_i$ in $W$.

We extend the length function $\ell:\waff\rightarrow\bn$ to a length
function $\ell:\wiff\rightarrow \bn$ by setting $\ell(\sigma w)=\ell(w)$
for $w\in \waff$ and $\sigma\in \Sigma$.%
%%%%%%%%%%%%%%%%%%%%%%%%%%%%%%%%%%%%%%%%%%%%%%%%%%
%%%%%%%%%%%%%%%%%%%%%%%%%%%%%%%%%%%%%%%%%%%%%%%%%%
%%%%%%%%%%%%%%%%%%%%%%%%%%%%%%%
%         Definiton of Demazure modules
%%%%%%%%%%%%%%%%%%%%%%%%%%%%%%%%%%%%%%%%%%%%%%%%%%
%%%%%%%%%%%%%%%%%%%%%%%%%%%%%%%%%%%%%%%%%%%%%%%%%%
%%%%%%%%%%%%%%%%%%%%%%%%%%%%%%%
\subsection{Definition of Demazure modules}\label{defndemazure}
For a dominant weight $\Lam\in \ph^+$ let $V(\Lam)$
be the (up to isomorphism) unique irreducible
highest weight module of highest weight $\Lam$.

Let $U(\Lbc)$ be the enveloping algebra of the
Iwahori subalgebra $\Lbc=\Lg \otimes t\bc[t] \oplus \Lb \otimes 1$,
and let $U(\Lhn)$ be the enveloping algebra of
$\Lhn=\Ln^+ \otimes \bc[t]\oplus \Lh \otimes t\bc[t] \oplus
\Ln^- \otimes t\bc[t]$.

Given an element
$w\in\waff/W_\Lam$, fix a generator $v_{w(\Lam)}$
of the line $V(\Lam)_{w(\Lam)}=\bc  v_{w(\Lam)}$
of $\Lhh$--eigenvectors in $V(\Lam)$ of weight $w(\Lam)$.
\begin{defn}
The $U(\Lhb)$--submodule $V_w(\Lam)=U(\Lhb)\cdot v_{w(\Lam)}$
generated by $ v_{w(\Lam)}$ is called the {\it Demazure submodule
of $V(\Lam)$} associated to $w$.
\end{defn}
\begin{rem}\rm
Since $v$ is an $\Lhh$--eigenvector,
we can also view the Demazure module $V_w(\Lam)$ as a cyclic  %\marginpar{ge\"andert}
$U(\Lbc)$--module or a cyclic $U(\Lhn)$--module generated by $ v_{w(\Lam)}$:
$$
V_w(\Lam)=U(\Lbc)\cdot v_{w(\Lam)}=U(\Lhn)\cdot v_{w(\Lam)}.
$$
\end{rem}
To associate more generally to every element $\sigma w\in\wiff=\Sigma\semi\waff$
a Demazure module, recall that elements in $\Sigma$ correspond to automorphisms
of the Dynkin diagram of $\Lhg$, and thus define an associated automorphism
of $\Lhg$, also denoted $\sigma$. For a module $V$ of $\Lhg$
let $V^\sigma$ be the module with the twisted action $g\circ v=\sigma^{-1}(g)v$.
Then for the irreducible module of highest weight $\Lam\in \ph^+$ we get
$V(\Lam)^\sigma=V(\sigma(\Lam))$.

So for $\sigma w\in \wiff=\Sigma\semi\waff$ we set
\begin{equation}\label{demazuredefn}
V_{w\sigma }(\Lam)=V_{w}(\sigma(\Lam))
\quad{\rm respectively}\quad
V_{\sigma w}(\Lam)=V_{\sigma w\sigma^{-1}}(\sigma(\Lam)).
\end{equation}
Recall that for a simple root $\al$ the Demazure module $V_{w\sigma }(\Lam)$
is stable for the associated subalgebra ${\mathfrak{sl}}_2(\al)$ if
and only if $s_\al w\sigma\le w\sigma\mod W_\Lam$ in the (extended)
Bruhat order. In  particular, $V_{w\sigma }(\Lam)$ is a $\Lg$--module
if and only if $s_i w\sigma\le w\sigma\mod W^{\rm aff}_\Lam$ for all $i=1,\ldots,n$.

We are mainly interested in Demazure modules
associated to the weight $\ell\Lam_0$ for $\ell\ge 1$.
In this case $W^{\rm aff}_\Lam=W$, so $\wiff/W=L$.
The Demazure module  $V_{t_{\nu(\mu^\vee)}}(\Lam_0)$
is a $\Lg$--module if and only if $\mu^\vee$ is an anti-dominant coweight,
or, in other words, $\mu^\vee=-\lam^\vee$ for some dominant coweight.
Since we will mainly work with these $\Lg$--stable Demazure modules,
to simplify the notation, we write in the following
%%%% %%%% %%%% %%%% %%%% %%%% %%%% %%%%
\begin{equation}\label{nummereins}
D(\ell, \lam^{\vee})\mbox{ for }\quad V_{t_{-\nu(\lam_{*}^\vee)}}(\ell\Lam_0)
\end{equation}
where $\lam_{*}^{\vee} = - w_0 (\lam^{\vee})$, the dual coweight of $\lam^{\vee}$.
This notation is justified by the fact that $D(\ell, \lam^{\vee})$ is,
considered as $\Lg$-module, far from being irreducible, but this $\Lg$-module
still has a unique
maximal highest weight: $\ell \nu(\lam^{\vee})$, i.e.,
if $V(\mu)$ is an irreducible $\Lg$-module of highest weight $\mu$
and $\Hom(V(\mu), D(\ell,\lam))\not=0$, then necessarily we have
$\ell\nu(\lam^\vee)-\mu$ is a non-negative sum of positive roots.
For more details on the $\Lg$-module structure of $D(\ell, \lam^{\vee})$
see also Theorem~\ref{fol-theorem} respectively \cite{FoL}.
%%%% %%%% %%%% %%%% %%%% %%%% %%%% %%%%
%%%%%%%%%%%%%%%%%%%%%%%%%%%%%%%%%%%%%%%%%%%%%%%%%%
%%%%%%%%%%%%%%%%%%%%%%%%%%%%%%%%%%%%%%%%%%%%%%%%%%
%%%%%%%%%%%%%%%%%%%%%%%%%%%%%%%
%          Demazure properties
%%%%%%%%%%%%%%%%%%%%%%%%%%%%%%%%%%%%%%%%%%%%%%%%%%
%%%%%%%%%%%%%%%%%%%%%%%%%%%%%%%%%%%%%%%%%%%%%%%%%%
%%%%%%%%%%%%%%%%%%%%%%%%%%%%%%%
\subsection{Properties of Demazure modules}
A description of Demazure modules
in terms of generators and relations
has been given by Joseph \cite{Joseph} (semisimple Lie algebras, characteristic zero)
and Polo \cite{Polo} (semisimple Lie algebras, characteristic free),
and Mathieu \cite{Mathieu} (symmetrizable Kac--Moody algebras).
We give here a reformulation for the affine case.
\begin{thm}[\cite{Mathieu}]
Let $\Lam\in \ph^+$ and let $w$ be an element of the affine
Weyl group of $\Lhg$. The Demazure module $V_w(\Lam)$ is as a $U(\Lhb)$-module
isomorphic to the following cyclic module, generated by $v \neq 0$ with the following
relations: for all positive roots $\beta$ of $\Lg$ we have
$$
\begin{array}{rcl}
(X_{\beta} \otimes t^{s})^{k_{\beta}+1}.v = 0 & \mbox{where } s\ge 0,
&k_{\beta}=max\{0,-\langle w(\Lam), (\beta+s\delta)^{\vee}\rangle\} \\
\\
(X_{\beta}^{-} \otimes t^{s})^{k_{\beta}+1}.v = 0 & \mbox{where } s>0, &
k_{\beta}=max\{0,-\langle w(\Lam),(-\beta + s\delta)^{\vee}\rangle\} \\
\\
(h \otimes t^{s}).v = 0 &\forall \, h \in \Lh,\ s>0,\\
\\
(h \otimes 1).v = w(\Lam)(h) v & \forall \, h \in \Lh,& d.v=w(\Lam)(d).v,\  K.v=\hbox{level}(\Lam) v
\end{array}
$$
\end{thm}
Let $\lam^{\vee}\in \chp^+$ be a dominant coweight.
We reformulate now the description of the Demazure modules above
for the Demazure modules $D(\ell, \lam^{\vee})$ we are interested in.
\begin{coro}\label{demazurecurrent}
As a module for the current algebra $\Lgc$,
$D(\ell, \lam^{\vee})$ is isomorphic to the cyclic $\Lgc$--module generated by a
vector $v$ subject to the following relations:
$$
\Ln^{+} \otimes \bc[t] . v = 0 \; , \; \Lh \otimes t\bc[t].v = 0 \; , \; h.v =\ell \nu(\lam^{\vee})(h)v
\mbox{ for all }h \in  \Lh,
$$
and for all positive roots $\beta\in\Phi^+$ one has
\begin{equation}\label{DEMAII}
(X_{\beta}^{-} \otimes t^s)^{k_{\beta}+1}.v = 0
\mbox{ where }
s\ge 0\mbox{ and } k_{\beta}= \ell \max\{0 ,-\langle \Lam_0+\nu(\lam^{\vee}),
(-\beta+s\delta)^{\vee}\rangle \}
\end{equation}
\end{coro}
\proof
Denote by $M$ the cyclic $U(\Lgc)$--module obtained by the relations above.
Recall (see (\ref{muaction})) that $t_{\nu(\lam^\vee)}(\ell\Lam_0)
= \ell\Lam_0+\ell \nu(\lam^\vee)$  and set  $\mu=t_{\nu(\lam^\vee)}(\ell\Lam_0)$.
Write $t_{\nu(\lam^\vee)}=w\sigma$ where $w\in \waff$ and $\sigma\in\Sigma$.
Set $\Lam=\sigma(\Lam_0)$, then the highest weight $\Lhg$-module
$V(\ell\Lam)$ has a unique line of $\Lhh$--eigenvectors
of weight $\mu$, let $v_\mu$ be a generator. Fix also a generator $v_{w_0(\mu)}$
of weight $w_0(\mu)$.
%So the Demazure module $V_{t_{\nu(\lam_*^\vee)}}(\ell\Lam_0)$ is generated
%as a $U(\Lhb)$--module by an extremal weight vector $v_{\lam^*}$ of weight
%$\ell\Lam_0+\ell \nu(\lam_*^\vee)$, and,
Restricted to the current algebra, we have $\mu\vert_\Lh=\ell \nu(\lam^\vee)$.
The submodule $U(\Lgc).v_{\mu}$ of $V(\ell\Lam)$ is the Demazure module
$D(\ell,\lam^\vee)$ because:
$$
D(\ell,\lam^\vee)= V_{t_{-\nu(\lam_{*}^\vee)}}(\ell\Lambda_0)= U(\Lhn)\cdot v_{w_0(\mu)}=U(\Lgc).v_{\mu}.
$$\noindent
Since $v_{\mu}$ is an extremal weight vector, using ${\mathfrak{sl}}_2$-representation theory
one verifies easily that $v_{\mu}$ satisfies the relations above. For example,
if the root is of the form $\beta+s\delta$, where $s\ge 0$ and $\beta\in\Phi^+$ is a positive root,
then the corresponding coroot is of the form $\beta^\vee+s' K$, $s'\ge 0$.
It follows that
$$
\langle \ell\Lam_0+\ell \nu(\lam^\vee), \beta^\vee+s' K\rangle =\ell s'+
\langle \ell \nu(\lam^\vee), \beta^\vee\rangle\ge 0,
$$
and hence $(\Ln^+\otimes\bc[t]) v_\mu=0$. So we have an obvious surjective
$\Lgc$--equivariant morphism $M\longrightarrow D(\ell,\lam^\vee)$, which
maps the cyclic generator $v$ to the cyclic generator $v_\mu$.

To prove that this map is an isomorphism it suffices to prove:
$\dim M\le\dim D(\ell,\lam^\vee)$.
The module $M$ is not trivial by the above, and the generator
$v\in M$ is a highest weight vector for the Lie subalgebra
$\Lg\subset \Lgc$. In fact, the relations imply that
the $\Lg$-submodule $U(\Lg).v\subseteq M$ is an irreducible, finite dimensional
highest weight $\Lg$-module $V(\ell\nu(\lam^\vee))\subseteq M$. So
we may replace for convenience the generator $v$ by a generator
$v'\in V(\ell\nu(\lam^\vee))$ of weight $w_0(\ell\nu(\lam^\vee))$,
i.e., we replace a $\Lg$-highest weight vector by a $\Lg$-lowest weight
vector. By construction, the following relations hold:
$$
\begin{array}{rrrl}
1)&(X_{\beta}^{+} \otimes t^s)^{k_{\beta}+1}.v' &=& 0\mbox{ where }
s\ge 0\mbox{; } k_{\beta} =\ell \max\{0 ,-\langle \Lam_0+ w_0(\nu(\lam^{\vee})),
(\beta+s\delta)^{\vee}\rangle \}\\
2)& (h\otimes 1).v'&=& \ell\nu(w_0(\lam^\vee))(h)v'
\mbox{ where } h\in\Lh.\\
3)&\Lh \otimes t\bc[t].v' &=& 0  \quad 4)\quad\Ln^- \otimes \bc[t].v' = 0 \\
\end{array}
$$
Now in $4)$  we have roots of the form $-\beta+s\delta$, where
$\beta$ is a positive root and  $s\geq0$. It follows
$$
\langle \Lam_0 + w_0(\nu(\lam^\vee)), -\beta^\vee+s' K\rangle = s'+
\langle - w_0(\nu(\lam^\vee)), \beta^\vee\rangle\ge 0,
$$
so we can reformulate $4)$ in the following way:
$$
4')\quad(X_{\beta}^{-} \otimes t^s)^{k_{\beta}+1}.v' = 0 \mbox{ where }
s\geq 0\mbox{ ; } k_{\beta}= \ell \max\{0 ,-\langle \Lam_0+w_0(\nu(\lam^{\vee})),
(-\beta+s\delta)^{\vee}\rangle \}
$$
Now $4)$ implies $M=U(\Lg\otimes\bc[t]).v'=U(\Lbc).v'$, and $1),2),3),4')$
show that the cyclic generator $v'$ for $M$ as $U(\Lbc)$--module
satisfies the same relations as the generator for the Demazure module
$D(\ell,\lam^\vee)$. Hence we have a surjective $U(\Lbc)$--module
homomorphism $D(\ell,\lam^\vee)\rightarrow M$, which finishes the proof.
\endpf
\begin{rem}\rm
We can easily extend the defining relations in Corollary \ref{demazurecurrent}
to an action of $\Lg\otimes\bc[t]\oplus\bc K$ by letting
$K$ act by $\ell$, the level of $\ell\Lambda_0$. This follows
immediately from (\ref{Lieklammer}) since in the current algebra
there are no elements of the form $x \otimes t^{-s}$, $s>0$.
\end{rem}
The $\Lg$--module structure of these special  Demazure modules
has been investigated in \cite{FoL}:
let $\lam^{\vee} = \lam^{\vee}_1 + \ldots + \lam^{\vee}_r$ be a sum of dominant
integral coweights for $\Lg$ and let $\ell\in \bn$.

\begin{thm}[\cite{FoL}]\label{fol-theorem}
As $\Lg$-modules the following are isomorphic
$$
D(\ell, \lam^{\vee}) \simeq D(\ell, \lam^{\vee}_1) \otimes \ldots \otimes D(\ell,\lam^{\vee}_r)
$$
\end{thm}
In this paper we will extend this isomorphism to an
isomorphism of $\Lgc$-modules by replacing the tensor product by the fusion product.
%%%%%%%%%%%%%%%%%%%%%%%%%%%%%%%%%%%%%%%%%%%%%%%%%%
%%%%%%%%%%%%%%%%%%%%%%%%%%%%%%%%%%%%%%%%%%%%%%%%%%
%%%%%%%%%%%%%%%%%%%%%%%%%%%%%%%
%          Weyl modules
%%%%%%%%%%%%%%%%%%%%%%%%%%%%%%%%%%%%%%%%%%%%%%%%%%
%%%%%%%%%%%%%%%%%%%%%%%%%%%%%%%%%%%%%%%%%%%%%%%%%%
%%%%%%%%%%%%%%%%%%%%%%%%%%%%%%%
\subsection{Weyl modules for the loop algebra}\label{weylmoduleloop}
The Weyl modules for the loop algebra
$\Lgl$ have been introduced in \cite{CP1}.
These modules are classified by $n$-tuples of
polynomials $\pi=(\pi_1,\ldots,\pi_n)$
with constant term 1, and they have the following universal
property: every finite dimensional cyclic $\Lgl$
highest weight module generated by a one-dimensional highest weight
space is a quotient of some $W(\pi)$ (for a more precise formulation see \cite{CP1}).
So these can be considered as maximal finite dimensional cyclic
representations in this class. 
A special class of tuples of polynomials is defined as follows:
fix $\lam=\sum_{j=1}^n m_j\om_j$ a dominant integral weight for $\Lg$
and a nonzero complex number $a\in \bc^*$, and set
\begin{equation}\label{specweyl}
\pi_{\lam,a}=((1-au)^{m_1},\ldots, (1-au)^{m_n}).
\end{equation}
The Weyl modules $W(\pi_{\lam,a})$ are of special interest because
\begin{enumerate}
\item it has been shown in \cite{CP1} that a Weyl module
$W(\pi)$ is isomorphic to a tensor product $\bigotimes_{j} W(\pi_{\lam_j,a_j})$
of Weyl modules corresponding to this special class of polynomials. 
\item the defining relations for the Weyl module $W(\pi_{\lam,a})$ reduce to 
(see \cite{CL}): $W(\pi_{\lam,a})$ is the cyclic module generated
by an element $w_{\lam,a}$, subject to the relations
$$
(\Ln^+\otimes\bc[t,t^{-1}])w_{\lam,a}=0,\ (h\otimes t^s)w_{\lam,a}=a^s \lam(h) w_{\lam,a},
\ (x_{\al_i}^-\otimes1 )^{m_i+1} w_{\lam,a}=0
$$
for all $h\in\Lh$, $1\le i\le n$, $s\in\bz$.
\end{enumerate}
In the following we denote by $\lam_\pi=\sum_i\deg\pi_i\om_i$ the weight
associated to a $n$-tuple of polynomials $\pi$. 
%It was conjectured in 
%\cite{CP1} that the dimension of a $W(\pi)$ depends only on $\lam_\pi$,
%more precisely: it is equal to $\prod_i (\dim W(\pi_{\om_i,1}))^{\deg\pi_i}$.
%As already pointed out in the introduction, it follows from the existence
%of the global crystal basis for the extremal weight modules \marginpar{Comments}
%\cite{Ka3} and the special structure of the crystal basis \cite{BN} that
%the dimension of the Weyl modules depend only
%on $\lam_{\pi_q}$ respectively $\lam_\pi$, and that 
%(whenever one has an appropriate lattice) the specialization of
%the quantum Weyl module is the classical Weyl module. These
%results were then used in \cite{CM3} to prove the 
%conjectured dimension formula for the Weyl modules
%for the quantum loop algebra, and, by specialization, hence
%also for the classical Weyl modules. 

\subsection{Weyl modules for the current algebra}
Let $\lam=\sum m_i \om_i$ be a dominant weight for $\Lg$.
A class of Weyl modules $W(\lam)$ has also been introduced
for the current algebra.  In terms of generators and relations one has:
\begin{defn}\label{weylcurrent}
Let $\lam$ be a dominant integral weight of $\Lg$, $\lam = \sum m_i \om_i$.
Denote by $W(\lam)$ the $\Lgc$-module generated by an element $v$ with the relations:
$$
\Ln^{+} \otimes \bc[t] .v = 0 \; , \;\Lh \otimes t\bc[t] .v = 0 \; , \; h.v = \lam(h).v  \; ,  \; (x_{\alpha_i}^{-}\otimes 1)^{m_i +1}.v = 0
$$
for all $h \in \Lh$ and all simple roots $\alpha_i$. This module is called
the Weyl module for $\Lgc$ associated to $\lam\in P^+$.
\end{defn}
The same proofs as those in \cite{CP1} show that $W(\lam)$ exists, is finite dimensional
and has the same universal property (see also \cite{CM1} and \cite{CL}).
\begin{rem}\label{allewurzeln}\rm
It follows easily that for all positive roots
$\beta$ the following relation holds in $W(\lam)$:
$$
(X_{\beta}^{- } \otimes 1)^{k_{\beta}+1}.v =0 \; \mbox{ for } \; k_{\beta} =  \lam(\beta^{\vee}).
$$
\end{rem}
For $a\in\bc^*$ consider the Lie algebra homomorphism
$\varphi_a$ defined as follows:
$$
\varphi_a:\Lgc\longrightarrow\Lgc,\quad x\otimes t^m\mapsto x\otimes (t+a)^m.
$$
Now $W(\pi_{\lam,a})$ is module for the loop algebra and hence
by restriction also a module for the current algebra. It has been shown 
in \cite{CP1,CM1}
that the twisted $\Lgc$-module
$$
{\varphi_a}^*(W(\pi_{\lam,a})),\mbox{ where the action is defined by } (x\otimes t^m)\circ_{\varphi_a} w
= (x\otimes (t-a)^m)w
$$
is a cyclic $\Lgc$-module satisfying the relations in Definition~\ref{weylcurrent}, so:
\begin{lem}
As a $\Lgc$-module, ${\varphi_a}^*W(\pi_{\lam,a})$ is a quotient of $W(\lam)$.
\end{lem}
In \cite{FKL} the so called higher level Weyl modules were introduced:
\begin{defn}
Let $W$ be a cyclic $\Lgc$-module, with fixed generator $w$.
We denote by $W^{[k]}$ the $\Lgc$ submodule of $W^{\otimes k}$ generated by $w^{\otimes k}$.

For a  dominant integral weight $\lam$ let $W(\lam)$ be
the Weyl module for the current algebra. The Weyl module of level $k$
corresponding to $\lam$ is defined as
$$
W(\lam)^{[k]}
$$
\end{defn}
\begin{rem}\label{higherdemazure}\rm
Let $V_{w}(\Lam)$ denote the Demazure submodule in the irreducible
highest weight $\Lhg$-module $V(\Lam)$ corresponding to the Weyl group element $w$.
Then
$$
V_{w}(\Lam)^{[k]} = V_{w}(k\Lam)
$$
\end{rem}
\begin{rem}\rm
\cite{FKL}
Let $V,W$ be cyclic $\Lgc$-modules, and suppose that $V$ is a quotient of $W$.
Then $V^{[k]}$ is a quotient of $W^{[k]}$.
\end{rem}
%%%%%%%%%%%%%%%%%%%%%%%%%%%%%%%%%%%%%%%%%%%%%%%%%%
%%%%%%%%%%%%%%%%%%%%%%%%%%%%%%%%%%%%%%%%%%%%%%%%%%
%%%%%%%%%%%%%%%%%%%%%%%%%%%%%%%
%          Fusion products
%%%%%%%%%%%%%%%%%%%%%%%%%%%%%%%%%%%%%%%%%%%%%%%%%%
%%%%%%%%%%%%%%%%%%%%%%%%%%%%%%%%%%%%%%%%%%%%%%%%%%
%%%%%%%%%%%%%%%%%%%%%%%%%%%%%%%
\subsection{Fusion products for the current algebra}
In this section we recall some facts on tensor products and fusion products of cyclic
$\Lgc$-modules.
Let $W$ be $\Lgc$-module and let $a$ be a complex number. 
Let $W_{a}$ be the $\Lgc$-module defined by the pullback 
$\varphi_{a}^*W$, so $x \otimes t^s$ acts as $x \otimes (t-a)^s$.
The following is well known:
\begin{lem}[\cite{FL1}]\label{cyclicshift} 
Let $W^1, \ldots, W^r$ be cyclic graded, 
finite-dimensional $\Lgc$-modules with cyclic vectors $w_1, \ldots, w_r$
 and let $C=\{c_1, \ldots ,c_r\}$ be pairwise distinct complex numbers.
Then $w_1\otimes \ldots \otimes w_r$ generates $W^1_{c_1} \otimes \ldots \otimes W^r_{c_r}$.
\end{lem}
The Lie algebra $\Lgc$ has a natural grading and an associated
natural filtration $F^{\bullet}(\Lgc)$, where $F^{s}(\Lg \otimes \bc[t])$ is defined
to be the subspace of $\Lg$-valued polynomials with degree smaller or equal $s$.
One has an induced filtration also on the enveloping algebra $U(\Lgc)$.
Let now $W$ be a cyclic module and let $w$ be a cyclic vector for $W$.
Denote by
$W_{s}$ the subspace spanned by the vectors of the form $g.w$,
where $g \in \, F^{s}(U(\Lgc))$, and denote the associated graded $\Lgc$--module
by $\hbox{\rm gr}(W)$
$$
\hbox{\rm gr}(W) = \bigoplus_{i\ge 0} W_s / W_{s-1}
\mbox{ where }W_{-1} = 0.
$$
As $\Lg$-modules, $W$ and $\hbox{\rm gr}(W)$ are naturally
isomorphic, but in general not as $\Lgc$-modules. 
%\marginpar{ Ist das wirklich so?}
\begin{defn}[\cite{FL1}]
Let $W^i$ and $ c_i$ as above in Lemma~\ref{cyclicshift} . 
The $\Lgc$-module
$$
W^1 \ast \ldots \ast W^r := \hbox{\rm gr}_{C}(W^1_{c_1} \otimes \ldots \otimes W^r_{c_r})
$$
is called the fusion product.
\end{defn}
\begin{rem}\rm  %\marginpar{ge\"andert}
It would be more appropriate to write $W^1_{c_1} \ast \ldots \ast W^r_{c_r}$
for the fusion product, since a priori the structure of the fusion product
depends on the choice of $C$. It has been conjectured in fact in \cite{FL1} that
the fusion product a) does not depend on the choice of the pairwise distinct 
complex numbers $C\in \bc^r $,
and b) is associative. This has been proved in the case $\Lg={\mathfrak{sl}}_n$
for various fusion products. In this paper we will prove the 
independence and the associativity property for the fusion product
of the Demazure modules $D(\ell,\lam)$, which justifies the fact that
we omit almost always the pairwise distinct complex numbers 
in the notation for the fusion product. In \cite{AK}
it is shown that the fusion product of Kirillov-Reshetikhin modules
of arbitrary levels is independent of the parameters.
\end{rem}
\begin{rem}\rm The case $r=1$ is of course not excluded. For example,
let $W$ be a graded cyclic $\Lgc$-module.  Let $C=\{c\}$, where $c\in \bc$,
then $\hbox{\rm gr}_{C}(W)\simeq W$ as $\Lgc$-modules.
\end{rem}
%%%%%%%%%%%%%%%%%%%%%%%%%%%%%%%%%%%%%%%%%%%%%%%%%
%%%%%%%%%%%%%%%%%%%%%%%%%%%%%%%%%%%%%%%%%%%%%
%          KR-modules
%%%%%%%%%%%%%%%%%%%%%%%%%%%%%%%%%%%%%%%%%%%%%%%%%
%%%%%%%%%%%%%%%%%%%%%%%%%%%%%%%%%%%%%%%%%%%%%
\subsection{Kirillov-Reshetikhin modules}\label{KRmodules}
In \cite{C2} (see also \cite{CM2}) for each multiple of a fundamental weight
$m \om_i$ a $\Lgc$-module $KR(m\om_i)$ has been defined.
These modules are called Kirillov-Reshetikhin module because in many cases
(Lie algebras of simply laced type or of classical type, see \cite{C2}) they can 
also be obtained from the quantum Kirillov-Reshetikhin module
by specialization and restriction to the current algebra.
\begin{defn}
Let $KR(m\om_i)$ be the $\Lgc$-module generated by a vector $v \neq 0$ with relations
\begin{equation}\label{KRI}
(\Ln^{+}\otimes \bc [t]).v =0 \; , \;(\Lh \otimes t\bc[t]).v =0 \; , \; hv = m\om_i(h) \; ,\;  h \in \Lh
\end{equation}
\begin{equation}\label{KRII}
(X^{-}_{\alpha_i})^{m + 1 }v =  (X^{-}_{\alpha_i}  \otimes t)v  = 0
\quad\hbox{\rm and\ }(X^{-}_{\alpha_j})v = 0 \; \mbox{for } j\neq i.
\end{equation}
\end{defn}
%%%%%%%%%%%%%%%%%%%%%%%%%%%%%%%%%%%%%%%%%%%%%%%%%%
%%%%%%%%%%%%%%%%%%%%%%%%%%%%%%%%%%%%%%%%%%%%%%%%%%
%%%%%%%%%%%%%%%%%%%%%%%%%%%%%%%
%         The quantum case
%%%%%%%%%%%%%%%%%%%%%%%%%%%%%%%%%%%%%%%%%%%%%%%%%% 
%%%%%%%%%%%%%%%%%%%%%%%%%%%%%%%%%%%%%%%%%%%%%%%%%%
%%%%%%%%%%%%%%%%%%%%%%%%%%%%%%%
\subsection{Quantum Weyl modules}\label{quantumweylmodules}
Let $U_{q}(\Lgl)$ be the
quantum loop algebra over $\bc(q)$, $q$ an
indeterminate, associated to $\Lg$ (see \cite{CP4}).
As in the classical case, one can associate finite-dimensional 
modules of $U_{q}(\Lgl)$ to $n$-tuples of
polynomials $\pi_{q}$ with constant term 1 and
coefficients in $\bc(q)$ (see \cite{CP1}). These modules are called
quantum Weyl modules. Again, the following universal property
holds: every highest weight module generated by a one-dimensional
highest weight space is a quotient of $W_{q}(\pi_q)$ for some $n$-tuple $\pi_q$
(see $\cite{CP1}$). Each such module has a unique irreducible quotient
which we denote by $V_{q}(\pi_q)$.

For such a $n$-tuple $\pi_{q}=(\pi_{q,1},\ldots,\pi_{q,n})$ set  
$\lam_{\pi_{q}}=\sum_i \deg \pi_{q,i}\om_i$, and let $\pi_{q,\om_i,1}$
be defined as in the classical case. 

The connection with Demazure modules is given by a theorem due to Kashiwara.
We state the theorem only for the simply laced type, but it holds
in much more generality.
\begin{thm}[\cite{Ka1}]\label{Kashiwarathm}
Let $\Lg$ be a simple Lie algebra of simply laced type, then
$$
\dim W_{q}(\pi_{q,\om_i,1}) = \dim D_{q}(1,\om_i^\vee)
$$
and $W_{q}(\pi_{q,\om_i,1})$ is irreducible.
\end{thm}
\begin{rem}\rm
The classical Demazure module $V_{w}(\Lam)$ (resp. $D(m,
\lam^{\vee})$) is the $q \rightarrow 1$ limit of the quantized
Demazure module $V_{q,w}(\Lam)$ (resp. $D_{q}(m,
\lam^{\vee})$).
\end{rem}

\begin{defn}
The $n$-tuple $\pi_q$ is called integral if all coefficients are
in $\ca$, and if the coefficient of the highest degree term
is in $\bc^*q^\bz$.
\end{defn}
Let $U_{\ca}(\Lgl)$ be
the $\ca$ subalgebra defined in $\cite{CP1}$.
It has been shown in \cite{CP1} that for an integral
$n$-tuple $\pi_q$ the corresponding quantum Weyl module
$W_{q}(\pi_q)$ admits a $U_{\ca}(\Lgl)$-stable
$\ca$-lattice $W_{\ca}(\pi_q)\subset W_{q}(\pi_q)$.

Further, let $\bc_1$ be the $\ca$-module with $q$ acting by $1$,
then $U(\Lgl)$ is a quotient of $U_{\ca}(\Lgl)\otimes \bc_1$,
and $\overline{W_q(\pi_q)}:=W_{\ca}(\pi_q)\otimes\bc_1$
becomes in a natural way a $U(\Lgl)$-module.

Let $\overline{\pi_q}$ be the $n$-tuple of
polynomials obtained by setting $q=1$, so the coefficients are in
$\bc$. 
The universality property of the Weyl modules implies \cite{CP1}:
\begin{lem}\label{spezialqgleicheins}
If $\pi_{q}$ is integral, then the $U(\Lgl)$-module
$\overline{W_q(\pi_q)}$ is a quotient of the classical Weyl module
$W(\overline{\pi_q})$.
\end{lem}
Actually, as outlined in the introduction, using results from global basis theory
one can show that one has equality in the lemma above. But
we need in the following only this weaker formulation.

It was already pointed out in \cite{CP2} that the cyclicity result of Kashiwara (\cite{Ka1}, Theorem 9.1,
see \cite{VV} for the simply laced case) for twisted tensor products of 
cyclic modules implies a lower bound for the dimension:
\begin{equation}\label{dimensioninequality}
\dim W_{q}(\pi_q)\ge \prod_i (\dim W_{q}(\pi_{q,\om_i,1}))^{\deg \pi_{q,i}}
\end{equation}
%%%%%%%%%%%%%%%%%%%%%%%%%%%%%%%%%%%%%%%%%%%%%%%%%%
%%%%%%%%%%%%%%%%%%%%%%%%%%%%%%%%%%%%%%%%%%%%%%%%%%
%%%%%%%%%%%%%%%%%%%%%%%%%%%%%%%
%          Connections between the modules
%%%%%%%%%%%%%%%%%%%%%%%%%%%%%%%%%%%%%%%%%%%%%%%%%%
%%%%%%%%%%%%%%%%%%%%%%%%%%%%%%%%%%%%%%%%%%%%%%%%%%
%%%%%%%%%%%%%%%%%%%%%%%%%%%%%%%
\section{Connections between the modules}\label{connections}
\subsection{Quotients}
We have some obvious maps between the Weyl modules for the current algebra
and certain Demazure modules.
\begin{lem}\label{demazurweylquotient}
Let $\lam^{\vee}$ be a dominant integral coweight of $\Lg$. For all $m\ge 1$,
the Demazure module $D(m, \lam^{\vee})$ is a quotient of the Weyl module
$W(m\nu(\lam^{\vee}))$.
\end{lem}
\proof
This follows immediately by comparing the relations for the
Weyl module in Definition~\ref{weylcurrent}
and the relations for the Demazure module in Corollary~\ref{demazurecurrent}
\endpf
\begin{lem}\label{fusion-lemma}
Let $\lam^{\vee}_i$, $i = 1, \ldots, r$ be dominant integral coweights,
let $\lam^{\vee} = \lam^{\vee}_1 + \ldots + \lam^{\vee}_r$,and let $a_1, \ldots, a_n$
be pairwise distinct complex numbers. Then
$$
D(1,\lam_1^{\vee})_{a_1} \ast \ldots \ast D(1,\lam_r^{\vee})_{a_r}
$$
is a quotient of $W(\nu(\lam^{\vee}))$.
\end{lem}
\proof  % \marginpar{ich glaub die Sterne mussten weg}
Let $v_i\in D(1,\lam_j^{\vee})$ be the cyclic generator as in Definition~\ref{weylcurrent}
and let $\nu(\lam_{i}^{\vee})=\sum_j m^i_j\omega_j$,
then the following relations hold:
$$
\Ln^{+} \otimes \bc[t].v_i  =0\;, \; (\Lh\otimes t\bc[t])v_i=0;\
h \otimes 1 .v_i = \nu(\lam_{i}^{\vee})(h)v_i \;,\; (x_{\alpha_j}^{-} \otimes 1)^{m^i_{j} + 1}.v_i = 0
$$
Let $\nu(\lam_{}^{\vee})=\sum_j m_j\omega_j$, then
the following relations
for the fusion product follow from the relations above:
$$
\Ln^{+} \otimes \bc[t].(v_i^{\otimes_{i=1}^r})  =0\;, \;
h \otimes 1 .(v_i^{\otimes_{i=1}^r})  = \nu(\lam_{}^{\vee})(h)(v_i^{\otimes_{i=1}^r })
\;,\; (x_{\alpha_j}^{-} \otimes 1)^{m^i_{j} + 1}.(v_i^{\otimes_{i=1}^r }) = 0.
$$
To see that all the relations of the Weyl module are satisfied in the fusion product,
it remains to show that $\Lh \otimes t\bc[t]$ annihilates $v_1 \otimes \ldots \otimes v_r$.
So (recall that $(h \otimes t^{k}).v_i = 0$ for $k> 0$) we have for $n\ge 1$:
$$
\begin{array}{rcl}
(h \otimes t^{n}).v_1 \otimes \ldots \otimes v_r & =& \sum\limits_{i} v_1 \otimes \ldots \otimes 
(h \otimes (t+c_j)^{n})v_i \otimes \ldots \otimes v_r\\%
& = & \sum\limits_{j} c_j^{n}\nu(\lam_{j}^{\vee})(h) v_1 \otimes \ldots \otimes v_r\\%
& = & (\sum\limits_{j} c_j^{n}\nu(\lam_{j}^{\vee})(h)) v_1 \otimes \ldots \otimes v_r\\
\end{array}
$$
By definition, this an element in the $n$-th part of the filtration,
but obviously the vector $v_1 \otimes \ldots \otimes v_r$ in also $0$-th part
of the filtration. Hence in the fusion product we have
$(h \otimes t^{n}).v_1 \otimes \ldots \otimes v_r=0$ for $n\ge 1$.
It follows: $(\Lh \otimes t\bc[t]).v_1 \otimes \ldots \otimes v_r = 0$,
which finishes the proof.
\endpf
%%%%%%%%%%%%%%%%%%%%%%%%%%%%%%%%%%%%%%%%%%%%%%%%%%%%%%%%
%%%%%%%%%%%%%%%%%%%%%%%%%%%%%%%%%%%%%%%%%%%%%%%%%%%%%%%%
%%%%%%%%%%%%%%KR-modules
%%%%%%%%%%%%%%%%%%%%%%%%%%%%%%%%%%%%%%%%%%%%%%%%%%%%%%%%
%%%%%%%%%%%%%%%%%%%%%%%%%%%%%%%%%%%%%%%%%%%%%%%%%%%%%%%%
\subsection{KR-modules}
\begin{thm}\label{KRgleichD}
For a fundamental coweight $\om_i^{\vee}$ let $d_i$ such that 
$d_i \om_i = \nu(\om_i^{\vee})$. The Kirillov-Reshetikhin module $KR(d_im\om_i)$
is, as $\Lgc$-module isomorphic to the Demazure module $D(m, \om_i^{\vee})$.
In particular, in the simply laced case (i.e., the root system is of type ${\tt A}_n, 
{\tt D}_n, {\tt E}_n$) all Kirillov-Reshetikhin modules are Demazure modules. 
\end{thm}
\begin{rem}\rm
The fact that $D(m, \om_i^{\vee})$ is a quotient of a  Kirillov-Reshetikhin module
has been already pointed out in \cite{CM2}. In the same paper Chari and Moura
have shown that $D(m, \om_i^{\vee})$ isomorphic to $KR(d_im\om_i)$ for
all classical groups using character calculations. Our proof
is independent of these results, and holds for all types, and gives an
alternative proof of the fact that these modules are finite dimensional.
\end{rem}
\proof
The fact that $D(m, \om_i^{\vee})$ is a quotient of $KR(d_im\om_i)$
is obvious by comparing the relations of the $KR$ modules with the 
relations of the Demazure module from Corollary~\ref{demazurecurrent}.

To show that the $KR$-modules above are quotients of Demazure modules,
it remains to verify that the relations~(\ref{KRII}) above imply the relations~(\ref{DEMAII})
in Corollary~\ref{demazurecurrent}. So let $\beta$ be a positive root, set $\al=\alpha_i$, $d=d_i$
and $\om=\om_i$. Note that $[h\otimes t^k,X_\beta^-\otimes t^\ell]
=\beta(h)X_\beta^-\otimes t^{\ell+k}$ 
implies: 
\begin{equation}\label{translaterelation}
(X_\beta^-\otimes t^s).v=0\quad\Rightarrow (X_\beta^-\otimes t^r).v=0\ \forall r\ge s.
\end{equation}
The fact that (\ref{DEMAII}) holds for the elements $X_\beta^-$ follows by ${\frak{sl}}_2$-theory. 
If $\langle \om,\bvee\rangle=0$, then (\ref{DEMAII}) holds
for all elements $X_\beta^-\otimes t^s,\ s\ge 0$, by $(\ref{translaterelation})$. 

Assume now $\langle \om,\bvee\rangle>0$
and consider an element of the form $X_\beta^-\otimes t^s$ for some $s\ge 1$. Let $\gamma\not=\al$
be a simple root, to verify the relation for $X_\beta^-\otimes t^s$ is equivalent
to verify it for $X_{s_\gamma(\beta)}^-\otimes t^s$. By replacing $\beta$ by ${s_\gamma(\beta)}$
if $\langle \beta,\gamma^\vee\rangle>0$, without loss of generality we may
assume that either $\beta=\alpha$, in which case the relations are satisfied, or 
$\beta\not=\alpha$ and $\alpha$ is the only simple root such that $\langle \beta,\avee\rangle>0$.

We have $\langle \beta,\avee\rangle=j$, $j=1,2,3$. So $\beta'=s_{\al}(\beta)=\beta-j\al$ and, if $t\ge j$,
then, up to a scalar,
$$
X_\beta^-\otimes t^s=[X_\al^-\otimes t,[...,[X_\al^-\otimes t,X_{\beta'}^-\otimes t^{s-j}]...].
$$
Except for the case where $\al,\beta'$ are two short roots in a root system of type ${\tt G}_2$,
the elements $X_\al^-\otimes t,X_{\beta'}^-\otimes t^{s-j}$ generate the nilpotent part of a Lie algebra
of type ${\tt A}_2$, ${\tt B}_2$ or ${\tt G}_2$. 

We consider first the case $\al,\beta'$ are short roots in a root system of type ${\tt G}_2$. Let 
$\gamma$ be the long simple root, then $\beta'=\gamma+\al$ and $\beta=\gamma+2\al=\om$. 
We have $X_\gamma^-v=0$, $(X_\al^-\otimes t)v=0$ and hence 
$(X_{\beta'}^-\otimes t)v=[X_\gamma^-, X_\al^-\otimes t]v=0$. In the same way one concludes
$(X_{\beta}^-\otimes t^2)v=0$.
Now using the commutation relations, one sees that $X_\beta.((X_\beta^-\otimes t)^k.v)=0$
for all $k\ge 0$. So if $(X_\beta^-\otimes t)^k.v\not=0$, then this a highest weight vector
for the Lie algebra generated by $X_\beta$ and $X_\beta^-$, and hence 
$(X_\beta^-\otimes t)^{3m+1}.v=0$. It follows that the elements $X_{\beta'}^-\otimes t^s,
X_{\beta}^-\otimes t^s$, $s\ge 0$, satisfy in this case the relations for the Demazure module
$D(m,\omega^\vee)$.

Suppose now $\al,\beta'$ form a basis of a root system of type ${\tt X}_2$, $X={\tt A,B,G}$.
Using the higher order Serre relations (see for example \cite{L}, Corollary 7.1.7), one sees
by induction that $(X_\al^-\otimes t).v=0$ implies for some constant $c\in\bc$:
$$
\begin{array}{rcl}
(X_\beta^-\otimes t^s)^m.v&=&([X_\al^-\otimes t,[...,[X_\al^-\otimes t,X_{\beta'}^-\otimes t^{s-j}]...])^m.v\\
&=&c\cdot (X_\al^-\otimes t)^{jm}(X_{\beta'}^-\otimes t^{s-j})^m.v
\end{array}
$$
Now if $X_{\beta'}^-\otimes t^{s-j}$ satisfies the relations for the Demazure module
in~(\ref{DEMAII}), then so does $X_\beta^-\otimes t^s$.

In the simply laced case this finishes the proof since the arguments above provide an
inductive method reducing the verification of the relations to the case either $\beta=\al$ or $s=0$,
and in both cases we know already that the relations hold. 
In the case $\Lg$ is of type ${\tt B}_n,{\tt C}_n$ or ${\tt F}_4$, the procedure
reduces the proof to the
cases 1) $\beta=\al$, 2) $s=0$ (now in these two cases the proof is finished), or
3) $\beta$ is a long root, $\alpha$ is a simple short root and
$\langle \beta,\avee\rangle=2$. In this case the relations have to be verified
for the root vector $X_\beta^-\otimes t$.

Now except for one case (in type ${\tt F}_4$) the pair $(\al, \beta'=\beta -2\al)$
is such that $\omega(\beta'^\vee)=0$, so $X_{\beta'}^-.v=0$ and hence 
$(X_{\beta}^-\otimes t^2).v=0$. Now as above, using the commutation relations
one sees that  $X_{\beta}(X_{\beta}^-\otimes t)^k.v=0$ for all $k$, so if
$(X_{\beta}^-\otimes t)^k.v\not=0$, then this is a highest weight vector
for the Lie subalgebra generated by $X_\beta^\pm$. It follows 
$(X_{\beta}^-\otimes t)^{m+1}.v=0$, and hence $X_{\beta}^-\otimes t$
satisfies the relations for the Demazure module. 

In the remaining case in type ${\tt F}_4$ (using reflections by simple roots
$\gamma\not=\al$) it suffices to consider
the pair where $\al=\al_3$ and $\beta=\al_1+2\al_2+4\al_3+2\al_4$ (notation
as in \cite{Bourbaki}). Now $\beta'$ is a positive root for which
we already know that the relations for the Demazure
module hold. (Using simple reflections $s_\gamma$, $\gamma\not=\al$,
to verify the relations for $\beta'$ is equivalent to verify them for $\al_2+2\al_3$,
this is a positive root of the type discussed above).
One has $\omega(\beta'^\vee)=2$,
so $(X_{\beta'}^-\otimes t^2).v=0$, $(X_{\beta'}^-\otimes t)^{m+1}.v=0$
and $(X_{\beta'}^-)^{2m+1}.v=0$, and the induction procedure shows
$(X_{\beta}^-\otimes t^4).v=0$, $(X_{\beta}^-\otimes t^3)^{m+1}.v=0$
and $(X_{\beta}^-\otimes t^2)^{2m+1}.v=0$.
It remains to show that $(X_{\beta}^-\otimes t)^{3m+1}.v=0$. Suppose
$u=(X_{\beta}^-\otimes t)^{3m+1}.v\not=0$, then, by ${\mathfrak{sl}}_2$-theory,
$X_\beta^{2m+2}u\not=0$. Now using the commutation relations, a simple induction
procedure shows that $X_\beta^{n}(X_{\beta}^-\otimes t)^{3m+1}.v$
is a linear combination of terms of the form
$$
(X_{\beta}^-\otimes t)^{3m+1-2n+j}(X_{\beta}^-\otimes t^2)^{n-2j}(X_{\beta}^-\otimes t^3)^jv.
$$
For $n=2m+2$ one has $j\ge m+2$, and hence all the terms vanish. It follows $u=0$.

Now for $\Lg$ of type ${\tt G}_2$ the induction procedure (respectively the arguments above
for the two short roots) reduce the proof of the relations to the cases of the root
vectors $X_\beta^-\otimes t$ and $X_\beta^-\otimes t^2$. Here
$\alpha$ is the short simple root, $\beta'$ is the long simple root and $\beta=s_\al(\beta')$.
Since $(X_\beta^-)^{3m+1}v=0$ and $(X_{\beta}^-\otimes t^3).v=0$,
the commutation relations as above show that 
$X_\beta^{m+2}(X_\beta^-\otimes t)^{2m+1}v=0$ and hence, by ${\mathfrak{sl}}_2$-theory,
$(X_\beta^-\otimes t)^{2m+1}v=0$.

Note that the root vectors $X_\beta^-$ and $X_{\beta-\al}^-$ commute. Since   
$(X_{\beta-\al}^-\otimes t^2)v=0$, we see that $X_\al(X_{\beta}^-\otimes t^2)^k.v=0$ for all
$k\ge 0$. So if $(X_{\beta}^-\otimes t^2)^k.v$ is nonzero, then this is a highest weight vector for the
Lie subalgebra generated by $X_\al^\pm$. Since $\langle 3m\om-k\beta,\al^\vee\rangle=3m-3k$,
it follows that $(X_{\beta}^-\otimes t^2)^{m+1}v=0$.
\endpf
%%%%%%%%%%%%%%%%%%%%%%%%%%%%%%%%%%%%%%%%%%%%%%%%%%%
%%%%%%%%%%%%%%%%%%%%%%%%
%          sl_2
%%%%%%%%%%%%%%%%%%%%%%%%%%%%%%%%%%%%%%%%%%%%%%%%%%
%%%%%%%%%%%%%%%%%%%%%%%%%%%%%%%%%%%%%%%%%%%%%%%%%%
%%%%%%%%%%%%%%%%%%%%%%%%%%%%%%%
\subsection{The $\mathfrak{sl}_{2}$-case}
Before we come to the proof of the main results,
let us recall the case $\Lg=\mathfrak{sl}_{2}$.
Note that for $\mathfrak{sl}_{2}$ we have $\nu(\lam^{\vee}) = \lam=\lam_*$.
Recall:
\begin{thm}[\cite{CP1}]
Let $\lam = m \om$, then dim $W(\lam)= 2^{m}$.
\end{thm}
As immediate consequence one obtains (already proved in \cite{CP1}, see also \cite{CL}):
\begin{thm}\label{sl_2_case}
For $\Lg = \mathfrak{sl}_{2}$ one has $W(\lam) \simeq D(1,\lam)$
as $\Lgc$-modules.
\end{thm}
\proof
The Demazure module is a quotient of the Weyl module, and by \cite{FoL}
and the theorem above one knows
that  $\dim D(1,\lam) = (\dim \,D(1,\om))^{m} = 2^{m}=\dim W(\lam)$.
\endpf
%%%%%%%%%%%%%%%%%%%%%%%%%%%%%%%%%%%%%%%%%%%%%%%%%%
%%%%%%%%%%%%%%%%%%%%%%%%%%%%%%%%%%%%%%%%%%%%%%%%%%
%%%%%%%%%%%%%%%%%%%%%%%%%%%%%%%
%          The simply laced case
%%%%%%%%%%%%%%%%%%%%%%%%%%%%%%%%%%%%%%%%%%%%%%%%%%
%%%%%%%%%%%%%%%%%%%%%%%%%%%%%%%%%%%%%%%%%%%%%%%%%%
%%%%%%%%%%%%%%%%%%%%%%%%%%%%%%%
\subsection{The simply-laced case}\label{simplylacedsection}
In this section let $\Lg$ be a simple simply laced Lie algebra, so $\Lg$ is of
type ${\tt A}_n,\,{\tt D}_n$ or ${\tt E}_n$.
Note that in this case $\nu(\lam^{\vee}) = \lam$.
We are now ready to prove
\begin{thm}\label{WeylgleichDemazure}
Let $\Lg$ be simply laced. Let $\lam^{\vee}\in \check{P}^+$ be a dominant integral coweight
for $\Lg$. The $\Lgc$-Weyl module $W(\lam)$ is isomorphic to the
Demazure module $D(1,\lam^{\vee})$.
\end{thm}
\begin{rem}\rm
In \cite{CL} the result has been proved for $\mathfrak{sl}_{n}$ by 
showing that the dimension conjecture of \cite{CP1} is true for the
classical Weyl module for $\Lg=\mathfrak{sl}_{n}$. Our approach 
is different and uses the relations defining a Demazure module.
On the other hand, we obtain a proof of the dimension conjecture of \cite{CP1} 
for the simply laced case by combining the result above 
with Theorem~\ref{fol-theorem}, see Proposition~\ref{dim-prop}.
\end{rem}
\proof
We know already that the Demazure module is a quotient of the Weyl module.
By comparing the defining relations in Corollary~\ref{demazurecurrent} and
in Definition~\ref{weylcurrent}, we see that to prove that this map is
an isomorphism, it is sufficient to show for
the Weyl modules that the following set of relations hold:
for all positive roots $\beta\in\Phi^+$ and all $s\ge 0$ one has
\begin{equation}\label{zusatzgleich}
(X_{\beta}^{-} \otimes t^s)^{k_{\beta}+1}.v = 0
\mbox{ where }
s\ge 0\mbox{ and } k_{\beta}=\max\{0 ,-\langle \Lam_0+\nu(\lam^{\vee}),
(-\beta+s\delta)^{\vee}\rangle \}
\end{equation}
Let $\beta$ be a positive root of $\Lg$, let $s \in \bn$ be a nonnegative integer
and set
$$
k = \mbox{max } \{0,\lam(\beta^{\vee})- s\}.
$$
Let $w_{\lam}\in W(\lam)$ be a generator of weight $\lam$. To prove (\ref{zusatzgleich}),
we have to show
$$
(X_{\beta}^{-}\otimes t^{s})^{k+1}.w_{\lam} = 0.
$$
Let ${\mathfrak{sl}}_{\beta}$ be the Lie subalgebra generated by $X_{\beta}^{-}, X_{\beta}, \beta^{\vee}$.
Let $V$ be the ${\mathfrak{sl}}_{\beta}\otimes \bc[t]$-submodule of $W(\lam)$
generated by $w_\lam$,
i.e., $V = U({\mathfrak{sl}}_{\beta}\otimes \bc[t]).w_{\lam}$.
Then $V$ satisfies obviously the
defining relations for the ${\mathfrak{sl}}_{\beta}\otimes \bc[t]$-Weyl module
$W_{\beta}(\lam(\beta^{\vee}))$ (see Remark~\ref{allewurzeln}),
so $V$ is a quotient of this Weyl module $W_{\beta}(\lam(\beta^{\vee}))$.
By Theorem \ref{sl_2_case} we know for the current algebra
${\mathfrak{sl}}_{\beta}\otimes \bc[t]$ that the Weyl module $W_{\beta}(\lam(\beta^{\vee}))$
is the same as the Demazure module $D_{\beta}(1,\lam(\beta^{\vee}))$. In particular, the defining
relations of $D_{\beta}(1,\lam(\beta^{\vee}))$ hold for the corresponding
generator of $W_{\beta}(\lam(\beta^{\vee}))$, and hence also for the corresponding
generator of $V$. It follows: $(X_{\beta}^{-}\otimes t^{s})^{k+1}.w_{\lam} = 0$.
\endpf
The following proposition is an immediate consequence of Theorem \ref{fol-theorem} and \ref{WeylgleichDemazure}.
\begin{prop}\label{dim-prop}
Let $\Lg$ be a simple, simply laced Lie algebra and
let $\lam^{\vee} = \sum m_i \om_i^{\vee}$ be a dominant integral coweight.
The dimension of $W(\lam)$ is
$$
\dim W(\lam) = \prod (\dim W(\om_i))^{m_i}
= \prod (\dim  D(1,\om_i^{\vee}))^{m_i}
$$
\end{prop}
We can now describe the current algebra module $\varphi^{*}_a(W(\pi_{\lam,a}))$
obtained as a pull back from the Weyl module for the loop algebra.
Here $\lam = \sum m_i \om_i$ and $\pi_{\lam,a}$ is the $n$-tuple of polynomials as
in section \ref{weylmoduleloop}.
\begin{prop}\label{pullback-prop}
Let $\lam$ be a dominant, integral weight for $\Lg$ of simply laced type. Then
$$
\varphi^{*}_a(W(\pi_{\lam,a})) \simeq W(\lam)
$$
\end{prop}
\proof
We know that $\varphi^{*}_a(W(\pi_{\lam,a}))$ is a quotient of $W(\lam)$, 
so it suffices to show that $\dim W(\pi_{\lam,a})\ge \dim W(\lam)$.
We have already seen that the specialization
$\overline{W_q(\pi_{q,\lam,a})}$ at $q=1$
of a quantum Weyl module is a quotient of
the Weyl module $W(\pi_{\lam,a})$ (see Lemma~\ref{spezialqgleicheins}).
By  Theorem~\ref{Kashiwarathm}, the
inequality~(\ref{dimensioninequality}) and
Proposition~\ref{dim-prop} it follows hence:
$$
\begin{array}{rcl}
\dim W(\pi_{\lam,a})\ge \dim W_q(\pi_{q,\lam,a}) &\ge&
\prod (\dim  W_{q}(\pi_{q,\om_i,1}))^{m_i} \\%=\prod (\dim  D_{q}(1,\om_i^\vee))^{m_i}\\
%&=&\prod (\dim  D_{q}(1,\om_i^\vee))^{m_i} \\
&=&\prod (\dim  D(1,\om_i^\vee))^{m_i}
=\prod (\dim  W(\om_i))^{m_i} =\dim W(\lam).
\end{array}
$$
\endpf
As an immediate and simple consequence we see:
\begin{coro}\label{specialcoro}
Let $\Lg$ be a simple Lie algebra of simply laced type, let $\lam$ be a dominant
weight (for $\Lg$), let $\pi$ (resp. $\pi_q$) be an $n$-tuple of polynomials in $\bc[u]$ (resp. in $\bc(q)[u]$) 
with constant term 1 such that $\lam=\lam_\pi=\lam_{\pi_q}$.
\begin{enumerate}
\item $\dim W(\lam) = \dim W(\pi) = \dim W_{q}(\pi_{q}) = \dim D(1,\lam^{\vee}) = \prod_i (\dim W(\om_i))^{m_i}$
\item If $\pi_{q}$ is integral, then 
$\overline{W_q(\pi_q)}\simeq W(\overline{\pi_q})$ as $U(\Lgl)$-modules.
\item The quantum Weyl module $W_q(\pi)$ is irreducible (note, the $\pi_i$ have complex
coefficients), and its specialization at $q=1$ is the Weyl module $W(\pi)$ for the classical loop algebra. 
\end{enumerate}
\end{coro}
\begin{rem}\rm
This proof of the dimension conjecture (in the simply laced case) can be seen 
as a more elementary alternative to the proof outlined in the introduction using
results from global basis theory.
\end{rem}
\proof
The first claim follows from Theorem~\ref{WeylgleichDemazure}, Proposition~\ref{dim-prop},
Proposition~\ref{pullback-prop}, the tensor product decomposition property 
(see section~\ref{weylmoduleloop}) and the specialization arguments outlined 
in \cite{CP2} (see section~\ref{quantumweylmodules}).

Now 2) is an immediate consequence of 1). To prove 3), let $m_{a^j}$ be the multiplicity of the 
root $a^j\in\bc^*$ of the polynomial $\pi_j(u)$.
The tensor product $W_q=\bigotimes_{j,a^j} W_q(\om_j,a^j)^{\otimes m_{a^j}}$ over all $j$ and 
all roots $a^j$ of $\pi_j(u)$ is irreducible by Theorem 9.2, \cite{Ka2}, it is again a highest weight 
module associated to the right $n$-tuple of polynomials, and has the right dimension by 1),
so it follows $W_q=W_q(\pi)$. The rest of the claim follows from 2).
\endpf
%\begin{coro}\label{quantsimple}
%%%%%%%%%%%%%%%%%%%%
We can now also prove the first step of Conjecture 1 in \cite{FoL} for
$\Lg$ of simply laced type. In the case of a multiple of a fundamental weight,
this provides a method to reconstruct the KR-module structure
for $U(\Lgl)$ from the $U(\Lgc)$-structure on the Demazure module.  
\begin{coro}
Let $D(m, \lam^{\vee})$ be a Demazure module of level $m$, corresponding to $\lam^{\vee}$.
Then $D(m, \lam^{\vee})$ can be equipped with the structure of a $U(\Lgl\oplus\bc K)$-module such that
the $\Lg$-module structure of $D(m, \lam^{\vee})$ coming from
the construction of the Demazure
module and the $\Lg$-module structure of $D(m, \lam^{\vee})$ obtained by the restriction of the
$U(\Lgl\oplus\bc K)$-module structure coincide.
\end{coro}
\proof 
As a $\Lgc$-module, $D(m, \lam^{\vee})$ is a quotient of $W(m \lam)$
(Lemma~\ref{demazurweylquotient}). Let $N(m\lam)$ be the kernel of the
map, so $D(m, \lam^{\vee}) \simeq W(m\lam)/N(m \lam)$
as $\Lgc$-modules. By Corollary~\ref{pullback-prop}, we know that 
$W(m\lam)$ is isomorphic to
$\varphi^{*}_1 W(\pi_{m\lam,1})$ as module for the current algebra.

Let $N_1(m\lam)=\varphi^{*}_{-1}N(m\lam)$ be the submodule 
of $W(\pi_{m\lam,1})$ corresponding to $N(m\lam)$.  Using \cite{CM1}, 
Proposition 3.3 (see also \cite{CP1}), one can show that $x \otimes t^{-s}$
operates as a linear combination of elements of  $U(\Lgc)$ on 
$W(\pi_{m\lam,1})$. So a $U(\Lgc)$-submodule of $W(\pi_{m\lam,1})$
is actually a $U(\Lgl)$-submodule. Since $K$ is
central (and operates trivially), we conclude that $N_1(m \lam)$ is a
$U(\Lgl \oplus \bc K)$-submodule of $W(\pi_{m\lam,1})$.

So the quotient $W(\pi_{m\lam,1})/N_{1}(m\lam)$ is a $U(\Lgl \oplus \bc K)$-module, isomorphic
to the Demazure module $D(m, \lam^{\vee})$ as vector space. Further, 
since $\varphi^*$ does not change the $\Lg$-structure of a $\Lgc$-module, 
we see that the $\Lg$-module structure
on $D(m, \lam^{\vee})$ and on the quotient $W(\pi_{m\lam,1})/N_{1}(m\lam)$ are identical.
\endpf
%\begin{rem}\rm
We conjecture that the corresponding statement also holds
for the quantum algebras and that the module admits 
a crystal basis as $U_q(\Lgl)$-module. Its crystal graph
should be obtained from the crystal graph of the quantum Demazure module just
by adding certain arrows with label zero. 

In the level 1 case we know
that we can identify $W_q(\pi_{q,\lam,1})$ with $D_q(1,\lam)$. To compare
the crystals, let $P_{cl}=P/\bz\delta$ be the quotient of the weight lattice by the
imaginary root and let $\psi:P\otimes_\bz \br\rightarrow P_{cl}\otimes_\bz \br$ 
be the projection of the associated real spaces. For a weight $\nu$ let
$\pi_\nu:[0,1]\rightarrow P\otimes_\bz \br$, $t\mapsto t\nu$, be the 
straight line path joining the origin with $\nu$, and let 
 $\psi(\pi_{\nu})$ be the image of the path in $P_{cl}\otimes_\bz \br$.
\begin{prop}
The crystal graph of $D_q(1,\lam)$ is obtained from the crystal
graph of $W_q(\pi_{q,\lam,1})$ by omitting certain arrows with label zero.
More precisely, let $B(\lam)_{cl}$ be the path model for
$W_q(\pi_{q,\lam,1})$ described in \cite{NS3}, then the crystal
graph of the Demazure module is isomorphic to the graph of the concatenation 
$\psi(\pi_{\Lambda_0})*B(\lam)_{cl}$.
\end{prop}
\proof
Write $t_{-\lam_{*}}$ as $w\sigma$, so $D_q(1,\lam)$ is the Demazure submodule
$V_{q,w}(\sigma(\Lambda_0))$. The path model theory  (see \cite{Li}) is independent 
of the choice of an initial path, we are going to 
choose an appropriate path. Instead of the straight line $\pi_{\sigma(\Lambda_0)}$ 
joining $0$ and $\sigma(\Lambda_0)$, consider 
the two straight line paths $\pi_{\Lambda_0}$ and $\pi_{-\lam_{*}}$ joining the origin
with $\Lambda_0$ respectively $-\lam_{*}$ in $P\otimes_\bz \br$. Let 
$\eta= \pi_{\Lambda_0}*\pi_{-\lam_{*}}$ be the concatenation of these two and
denote by $B(\eta)$ the set of paths generated by applying the root operators to $\eta$.
By \cite{Li}, section 6, $\eta$ is linked for arbitrary $L$ to the straight line path 
$\pi_{\Lambda_0-\lam_{*}}$, which is an LS-path of shape  $\sigma(\Lambda_0)$. 
It follows that the two path models are isomorphic, and hence: a) the crystal associated to
the set of paths $B(\eta)$ is isomorphic to the crystal of $V_{q}(\sigma(\Lambda_0))$,
and b) in $B(\eta)$ there exists a unique path $\pi_0$ contained in the dominant 
Weyl chamber and ending in $\sigma(\Lambda_0)$. Denote by $B(\eta)_{cl}$ the image
of this set of paths under the projection $\psi$.  The root operators $e_\alpha,f_\alpha$,
$\alpha$ a simple root for $\Lhg$, are still well defined on paths in $P_{cl}\otimes_\bz \br$
since $\delta$ vanishes on all coroots. In fact, the operators
commutes with the map $\psi$. So the uniqueness of $\eta$ (as path contained in the dominant 
Weyl chamber) implies that $\psi$ induces a bijection between the crystals
$B(\eta)$ and $B(\eta)_{cl}$.

Let $B(\lam)$ be the set of all LS-paths of shape $\lam$ and denote by $B(\lam)_{cl}$
the image of this set under the projection $\psi$. Combining part 3 of the theorem above 
with the result of Naito and Sagaki in \cite{NS3}, we see that $B(\lam)_{cl}$ is a 
combinatorial model for the crystal of the Weyl module $W_q(\pi_{q,\lam,1})$. 
The concatenation $\psi(\pi_{\Lambda_0})*B(\lam)_{cl}$ in $P_{cl}\otimes_\bz \br$ 
provides a set of paths stable under all root operators $e_\alpha$,
$\alpha$ a simple root, and $f_{\alpha_i}$, $i=1,\ldots,n$. 

To describe in $B(\eta)_{cl}$ the set of path corresponding to the Demazure module
$D(1,\lam)$, recall that the latter is the union of all paths of
the form $f_{\al_{j_1}}^{n_1}\cdots f_{\al_{j_t}}^{n_t}\pi_0$, where 
$t_{-\lam_{*}}=w\sigma$. $w=s_{\al_{j_1}}\cdots s_{\al_{j_t}}$ is a 
reduced decomposition and $n_i\in\bn$. Recall that $\pi_0$ is of 
the form $\pi_{\Lambda_0}*\pi'$, where $\pi'\in B(\lam)$. Since we work modulo 
$\delta$ and $\lam$ is a level zero weight, $e_\al\pi_0=0$ for all simple
roots of $\Lhg$ means that $\pi'$ is a path completely contained in the 
fundamental alcove of the root system of $\Lg$. This path is obtained 
from the straight line path $\pi_\lam$ by folding it successively back into
the alcove (in the same way as in \cite{Li1}, proof of the PRV-conjecture).
Next consider the sequence of turning points. 

If $\lam$ is regular and generic (i.e., $\lam\not=m\mu$ for some $m\ge 2$, 
$\mu\in P^+$), then this are exactly the points where $\pi_0$ meets the 
codimension one faces of the fundamental alcove $\Delta_f$, and the 
corresponding product of the simple reflections 
is exactly a reduced decomposition of $w'$, where $t_{\lam}=w'\sigma$.
We get a reduced decomposition of $w$ (recall, $t_{-\lam_{*}}=w\sigma$)
by multiplying the given reduced decompositions with appropriate simple
reflection $s_{\al_i}$, $i\ge 1$. By the choice of the reduced decomposition
above, the paths $f_{\al_{j_1}}^{n_1}\cdots f_{\al_{j_t}}^{n_t}\pi_0$
are all of the form $\psi(\pi_{\Lambda_0})*\pi'$, where $\pi'\in B(\lam)$.

The same holds also in the general case, only that the turning points 
are not anymore associated to just one simple reflection an 
element of maximal length in a coset $W'/W''$ of subgroups of $\waff$.
Here $W', W''$ are associated to the turning point and the path $\pi_0$,
for details in terms of galleries see for example \cite{GL}, Example 4. 

So the set of paths in the path crystal of $V_{q}(\sigma(\Lambda_0))$ corresponding
to the subcrystal of $D(1,\lam)$ is a subset of $\psi(\pi_{\Lambda_0})*B(\lam)_{cl}$.
By the equality of the number of elements, the two sets have in fact to be equal.
\endpf 
%%%%%%%%%%%%%%%%%%%%%%%%%%%%%%%%%%%%%%%%%%%%%%%%%%
%%%%%%%%%%%%%%%%%%%%%%%%%%%%%%%%%%%%%%%%%%%%%%%%%%
%%%%%%%%%%%%%%%%%%%%%%%%%%%%%%%??????????????????????????????
%          Demazure and Fusion module
%%%%%%%%%%%%%%%%%%%%%%%%%%%%%%%%%%%%%%%%%%%%%%%%%%
%%%%%%%%%%%%%%%%%%%%%%%%%%%%%%%%%%%%%%%%%%%%%%%%%%
%%%%%%%%%%%%%%%%%%%%%%%%%%%%%%%
\subsection{Demazure modules as fusion modules}
In this section let $\Lg$ be a simple Lie algebra of arbitrary type. So, unless 
it is explicitly mentioned, in this section we do not assume that $\Lg$ is 
necessarily simply laced.
\begin{thm}\label{demazure-zerlegung}
Let $\lam^{\vee} = \sum_{i=1}^s \lam_i^{\vee}$ be a sum of dominant integral coweights
and let $c_1, \ldots, c_s$ be pairwise distinct complex numbers, then
$$
D(1,\lam^{\vee})\simeq D(1,\lam_1^{\vee}) \ast  \ldots \ast D(1,\lam_s^{\vee}) 
$$
as modules for the current algebras $\Lg \otimes \bc[t]$.
\end{thm}
In the simply laced case we have of course equivalently:
\begin{coro}\label{weylfusionsimplylaced}
Let $\Lg$ be a simple simply laced Lie algebra.
For $\lam = \sum\limits_{i=1}^{s} \lam_i$, $\lam_i$ dominant integral coweights, 
and $c_1,  \ldots , c_s$ pairwise distinct complex numbers:
$$
W(\lam) \simeq W(\lam_1) \ast \ldots \ast W(\lam_s)
$$
as $\Lgc$-modules
\end{coro}
\begin{rem} \rm
For $\Lg={\mathfrak{sl}}_n$ and the $\lam_i$, $i=1,\ldots, s$, all
fundamental weights, the theorem above (and its corollary) was proved 
by Chari and Loktev in \cite{CL}.
\end{rem}
\begin{coro}\label{demazure-zerlegunghigh}
Let $\Lg$ be again a simple Lie algebra of arbitrary type and
let $\lam^{\vee} = \sum_{i=1}^s \lam_i^{\vee}$ be a sum of dominant integral coweights
and let $c_1, \ldots, c_s$ be pairwise distinct complex numbers, then for all $k\ge 1$
$$
D(k,\lam^{\vee})\simeq D(k,\lam_1^{\vee}) \ast  \ldots \ast D(k,\lam_s^{\vee})
$$
as modules for the current algebras $\Lg \otimes \bc[t]$.
\end{coro}
As obvious consequences we have:
\begin{coro}
\begin{enumerate}
\item The fusion product of the Demazure modules $D(k,\lam_j)$ is associative
and independent of the choice of the pairwise distinct complex numbers
$\{c_1, \ldots, c_s\}$.
\item  Let $d_i$ be as in Theorem~\ref{KRgleichD}, then
$KR(md_i\om_i)$ is the $m$-fold fusion product $KR(d\om_i)^{*m}$.
\end{enumerate}
\end{coro}
\noindent{\it Proof of Corollary~\ref{demazure-zerlegunghigh}.\/}
It follows from Remark~\ref{higherdemazure} and Theorem~\ref{demazure-zerlegung}
that 
$$
D(k,\lam^{\vee})=D(1,\lam^{\vee})^{[k]}\simeq
(D(1,\lam_1^{\vee}) \ast  \ldots \ast D(1,\lam_s^{\vee}))^{[k]}.
$$
By Proposition~2.10 in \cite{FKL} the latter is a quotient of 
$D(1,\lam_1^{\vee})^{[k]} \ast  \ldots \ast D(k,\lam_s^{\vee})^{[k]}
= D(k,\lam_1^{\vee}) \ast  \ldots \ast D(k,\lam_s^{\vee})$.
The dimension formula (Theorem~\ref{fol-theorem}) 
implies again that the map is an isomorphism.
\endpf
\noindent{\it Proof of Theorem~\ref{demazure-zerlegung}.\/}
In the simply laced case the result follows immediately from the
equality of Demazure and Weyl modules:
the right hand side is a Weyl module by Theorem~\ref{WeylgleichDemazure},
and the left hand side is a quotient of this Weyl module by
Lemma \ref{fusion-lemma}. Now by Theorem \ref{fol-theorem} the
dimension of both modules is equal, which finishes the proof.

In the general case we need to use the defining equations for Demazure module
(see Corollary~\ref{demazurecurrent}). In the proof of Lemma~\ref{fusion-lemma}
we have already seen that the fusion module:
$$
D(1,\lam_1^{\vee}) \ast \ldots \ast D(1,\lam_r^{\vee})
$$
is a quotient of the Weyl module and hence satisfies the relations:
$$
\Ln^{+} \otimes \bc[t].(v_i^{\otimes_{i=1}^r})  =0\;, \;
h \otimes 1 .(v_i^{\otimes_{i=1}^r})  = \nu(\lam_{}^{\vee})(h)(v_i^{\otimes_{i=1}^r })
\hbox{\rm\ and\ }\Lh \otimes t\bc[t](v_i^{\otimes_{i=1}^r }) = 0.
$$
Let now $\beta\in\Phi^+$ be a positive root. The following lemma
implies that the fusion product is a quotient of the Demazure module. Since
both have the same dimension by Theorem~\ref{fol-theorem}, they are 
isomorphic, which finishes the proof.
\endpf
\begin{lem}\label{eqdemazuregeberal}
\begin{equation}\label{demazuregeneral}
(X_\beta^-\otimes t^s)^{k_\beta+1}(v_i^{\otimes_{i=1}^r})=0\mbox{\ for\ }
k_\beta=\max\{ 0,\langle\Lambda_0+\nu(\lambda^\vee),(-\beta+s\delta)^\vee\rangle\},
\end{equation}
\end{lem}
The proof of Lemma~\ref{eqdemazuregeberal} is by reduction to the
$\widehat{\mathfrak{sl}}_2$--case. Note that in this case we know
already that Theorem~\ref{demazure-zerlegung} and 
Corollary~\ref{demazure-zerlegunghigh} hold.

We fix first some notation.
For a positive root $\beta\in\Phi^+$ 
let $Z_\beta\subset \Lhg$ be the Lie subalgebra generated 
by the root spaces $\Lhg_{\pm \beta + s\delta}$, $s\in\bz$,
the elements in the Cartan subalgebra $(\pm \beta \pm s\delta)^\vee$,
and the derivation $d$. Then $Z_\beta$ is an affine Kac-Moody algebra
isomorphic to $\widehat{\mathfrak{sl}}_2$ with Cartan subalgebra 
$\Lhh_\beta=\langle\beta^\vee, \eps K,d\rangle_\bc$, 
where $\eps=(\beta^\vee,\beta^\vee)/2$ (see 
equation (\ref{coroots}))
is $1$ if $\beta$ and $\Theta$ have the same length, and
$\eps=2$ or $3$ if $\beta$ is a short root. 
Set 
$$
\Lhn_\beta^\pm=\Lhn^\pm\cap Z_\beta\quad\hbox{\rm and}\quad
\mathfrak{sl}_2(\beta)\otimes\bc[t]=Z_\beta\cap \Lg\otimes\bc [t]
$$
Write $t_{\nu(\lam^\vee)}=w\sigma$ 
(see equations~(\ref{demazuredefn}) and (\ref{nummereins})),
so $w\in\waff$. Set $\mu=w(\sigma(\Lambda_0))$ and let
$v_\mu\in V(\sigma(\Lambda_0))$ be an extremal weight vector
of weight $\mu$. 
The submodule $M=U(Z_\beta)v_\mu\subset V(\sigma(\Lambda_0))$
is an irreducible (since $v_\mu$ is an extremal weight vector) 
$Z_\beta$-submodule, say $M=V^\beta(\Omega )$ is the
$Z_\beta$-representation of highest weight $\Omega$. The subspace
\begin{equation}\label{betademazure}
M(\nu(\lam^\vee)):=U(\mathfrak{sl}_2(\beta)\otimes\bc[t]).v_\mu
=U(\Lhn_\beta^+)s_\beta(v _\mu)
\end{equation}
is then a Demazure module, stable under $U(\mathfrak{sl}_2(\beta)\otimes\bc[t])$.
Now $V(\sigma(\Lambda_0))$ is a level one module for $\Lhg$,
but the irreducible $Z_\beta$-submodule $V_\beta$ is a 
level $\eps$-module for the affine Kac-Moody algebra 
$Z_\beta\simeq \widehat{\mathfrak{sl}}_2$
(recall, the canonical central element of $Z_\beta$ is $\eps K$).
We need the following more precise statement:
\begin{lem}\label{hilfslevellemma}
As ${\mathfrak{sl}}_2(\beta)\otimes \bc[t]$-module, the submodule 
$M(\nu(\lam^\vee))$ is isomorphic to
$D(\eps,m\om^\beta)$, where $m=\nu(\lam^\vee)(\beta^\vee)/\eps$ and $\om^\beta$ denotes
the fundamental weight for the Lie algebra ${\mathfrak{sl}}_2(\beta)$. 
\end{lem}
\par\noindent{\it Proof of Lemma~\ref{hilfslevellemma}.\/}
The first step is to show that the highest weight $\Omega$ is a multiple of a fundamental
weight for $Z_\beta$. The only non-trivial case is when $\Theta$ and $\beta$
have different lengths. We show first that in this case:
\begin{equation}\label{epsmodularitaet}
\nu(\lam^\vee)(\beta^\vee)\equiv 0\bmod \eps\bz.
\end{equation}
To prove this, recall that for $\lam^\vee=\sum_i m_i\om_i^\vee$ one has
$\nu(\lam^\vee)=\sum_i m_i \frac{a_i}{a_i^\vee}\om_i$. If $\al_i$
is a short root, then $\frac{a_i}{a_i^\vee}=\eps$, so 
$\frac{a_i}{a_i^\vee}\om_i(\beta^\vee)\equiv 0\bmod \eps\bz$.
Now a case by case consideration shows that if $\al_i$ is a long simple
root and $\beta$ is a short positive root, then again 
$\frac{a_i}{a_i^\vee}\om_i(\beta^\vee)\equiv 0\bmod \eps\bz$. 

Now $s_\beta(\nu(\lam^\vee)\vert_{\Lhh_\beta})\equiv t_{\eta}(\Omega)\bmod\bz\Lam_0$ 
( respectively $t_{\eta}(\sigma(\Omega))\bmod\bz\Lam_0$)
for some ${\mathfrak{sl}}_2(\beta)$-weight $\eta$. Since $t_{\eta}(\Omega)=\Omega+\epsilon\eta$
respectively  $t_{\eta}(\sigma(\Omega))=\sigma(\Omega)+\epsilon\eta$, 
it follows
\begin{equation}%\label{epsmodularitaet}
\Omega(\beta^\vee)\equiv 0\bmod \eps\bz\quad
\hbox{\rm respectively}\quad \sigma(\Omega)(\beta^\vee)\equiv 0\bmod \eps\bz.
\end{equation}
But this is only possible if $\Omega=\eps \Lam^\beta_0$ or $\Omega=\eps \Lam^\beta_1$
as highest weight for the irreducible
$Z_\beta\simeq \widehat{\mathfrak{sl}}_2$--representation $M$, and
hence $M(\nu(\lam^\vee))\simeq D(\eps,m\om^\beta)$ for some $m$. 
Since $\nu(\lam^\vee)(\beta^\vee) =(\eps\Lam^\beta_0+\eps m\om^\beta)(\beta^\vee)$, 
it follows that $m=\nu(\lam^\vee)(\beta^\vee)$.
\endpf
\par\noindent{\it Proof of Lemma~\ref{eqdemazuregeberal}.\/}
For each of the Demazure modules $D(1,\lam_i^\vee)$ denote by
$M(\nu(\lam_i^\vee))$ the $Z_\beta$--Demazure submodule generated by $v_i$,
as in (\ref{betademazure}). By the lemma above we have 
$M(\nu(\lam_i^\vee))\simeq D(\eps,m_i\om^\beta)$,
where $m_i=\nu(\lam_i^\vee)(\beta^\vee)/\eps$. Taking the tensor product, we get
an embedding
$$
M(\nu(\lam_1^\vee))\otimes\cdots\otimes M(\nu(\lam_s^\vee))
\hookrightarrow 
D(1,\lam_1^\vee)\otimes\cdots\otimes D(1,\lam_s^\vee)
$$
Now the filtration on $M(\nu(\lam_1^\vee))\otimes\cdots\otimes M(\nu(\lam_s^\vee))$
as ${\mathfrak{sl}}_2(\beta)\otimes \bc[t]$--module is compatible with the
filtration of $D(1,\lam_1^\vee)\otimes\cdots\otimes D(1,\lam_s^\vee)$ as
$\Lg\otimes\bc[t]$--module, so we get an induced map
$$
M(\nu(\lam_1^\vee))*\cdots * M(\nu(\lam_s^\vee))
\longrightarrow 
D(1,\lam_1^\vee)*\cdots * D(1,\lam_s^\vee)
$$
Since we are in the simply laced case, we know 
by Corollary~\ref{demazure-zerlegunghigh} that the left
$\mathfrak{sl}_2(\beta)\otimes\bc[t])$-module is isomorphic
to $M(\nu(\lam^\vee))$. By construction, the
generator of this module satisfies the equation~(\ref{demazuregeneral}),
and hence also the image $v_i^{\otimes_{i=1}^{s}}$ satisfies 
equation~(\ref{demazuregeneral}).
\begin{rem}\rm
It is not anymore true that $W(\nu(\lam^{\vee})) \simeq D(1, \lam^{\vee})$.
As a counter example consider $\Lg$ of type ${\tt C}_2$ and
take $\lam^{\vee} = \om_1^{\vee}$. Note that $\nu(\om_1^\vee)=2\om_1$. 
By \cite{FoL}, $D(1, \om_1^{\vee})$ has dimension $11$,
and by \cite{C2}, $KR(\om_1)$ has dimension $4$. The fusion product
$KR(\om_1)\ast KR(\om_1)$ is a quotient of $W(2\om_1)$, so 
$\dim W(2\om_1)\geq 16 > \dim D(1, \om_1^{\vee})$.
\end{rem}
We conjecture that Corollary~\ref{weylfusionsimplylaced} also holds in the
non-simply laced case:
\begin{conj}
Let $\lam = \sum \lam_i$ be a sum of dominant integral weights, 
$c_1, \ldots, c_n$ be pairwise distinct complex numbers, then
$$
W(\lam) \simeq W(\lam_1) \ast  \ldots \ast W(\lam_n)
$$
\end{conj}
%%%%%%%%%%%%%%%%%%%%%%%%%%%%%%%%%%%%%%%%%%%%%%%%%%
%%%%%%%%%%%%%%%%%%%%%%%%%%%%%%%%%%%%%%%%%%%%%%%%%%
%%%%%%%%%%%%%%%%%%%%%%%%%%%%%%%
%          Limit constructions
%%%%%%%%%%%%%%%%%%%%%%%%%%%%%%%%%%%%%%%%%%%%%%%%%%
%%%%%%%%%%%%%%%%%%%%%%%%%%%%%%%%%%%%%%%%%%%%%%%%%%
%%%%%%%%%%%%%%%%%%%%%%%%%%%%%%%
\section{Limit constructions}\label{seclimitconstruction}
In this section we start with a simple Lie algebra $\Lg$ of arbitrary type.
We want to reconstruct the $\Lgc$-module structure of the
irreducible highest weight $U(\Lhg)$-module $V(l\Lam_0)$
as a direct limit of fusion products of Demazure modules.

In \cite{FoL} we have given such a construction of
the $\Lg$-module structure of $V(l\Lam_0)$ as
a semi-infinite tensor product of finite dimensional
$\Lg$-module. In this section we want to extend this
construction 
%\marginpar{ge\"andert}
to the $U(\Lgc)$-module structure by
replacing the tensor product by the fusion product.

We need first a few facts about inclusions
of Demazure modules. Set $\tilde{\Lb} = \Lh \oplus \Lhn^{+} \oplus \bc K$,
and, as before, we denote by $\waff $ the affine Weyl group.
Let $\Lam$ be an integral dominant weight for $\Lhg$.
We fix for all $w\in\waff /\waff _\Lam$ a generator $v_w$
of the line of weight $w(\Lam)\subset V(\Lam)$. Denote
$V_w(\Lam)=U(\Lhb).v_w$ the Demazure module and
let $\iota_w:V_w(\Lam)\hookrightarrow V(\Lam)$ the
inclusion.
\begin{lem}\label{unique-embedding}
Let $\Lam$ be an integral dominant weight for $\Lhg$.
Given $w\in\waff /\waff _\Lam$, there is a unique
(up to scalar multiplication) nontrivial morphism of $U(\tilde{\Lb})$-modules
$$
V_{w}(\Lam) \longrightarrow V(\Lam).
$$
In fact, this morphism is, up to scalar multiples, the
canonical embedding of the Demazure module.
\end{lem}
\proof
We want to prove that, up to scalar multiples,
$\iota_w:V_{w}(\Lam) \longrightarrow V(\Lam)$
is the only nontrivial morphism of $U(\tilde{\Lb})$-modules.
The proof is by induction on the length of $w$.

For $w = id$, the Demazure module is one-dimensional.
The generator $v$ is killed by $U(\Lhn^{+})$, so its image in
$V(\Lam)$ is a highest weight vector.
But such a vector is unique (up to scalar multiple) in $V(\Lam)$, and
hence there exists, up to a scalar multiples, only one such
morphism.

Suppose now $\ell(w)\ge 1$, and let $\tau = s_{\alpha} w$, $\alpha$ a simple root,
be such that  $\tau <  w$, and let $\varphi : V_{w}(\Lam) \longrightarrow V(\Lam)$
be a non trivial $U(\tilde{\Lb})$-equivariant morphism.  Let $v_w$ be a generator
of the weight space in $V_{w}(\Lam)$ corresponding to the weight $w(\Lam)$,
and set $m_{\alpha} = w(\Lam)(\alpha^{\vee})$.
Then $(x_{\alpha})^{m_{\alpha}}.v_{w} \neq 0$, but
$(x_{\alpha})^{m_{\alpha}+1}.v_{w} = 0$.

Now $\varphi$ is an $U(\tilde{\Lb})$-morphism, so the image
$\varphi(v_{w})\in V(\Lam)$ is again an eigenvector for $\Lh\oplus \bc K$
of weight $w(\Lam)\vert_{\Lh\oplus\bc K}$. Since $V(\Lam)$ is a $\Lhg$-module,
${\mathfrak{sl}}_2$-representation theory implies
$(x_{\alpha})^{m_{\alpha}}.\varphi(v_{w}) \neq 0$, and since
 $\varphi$ is an $U(\tilde{\Lb})$-morphism, we have
$(x_{\alpha})^{m_{\alpha}+1}.\varphi(v_{w}) = 0$.

Now $(x_{\alpha})^{m_{\alpha}}.v_{w}$ is a generator of the Demazure module 
$V_{\tau}(\Lam)\subset V_w(\Lam)$, so $\varphi\vert_ {V_{\tau}(\Lam)}$
provides a non-trivial $U(\tilde{\Lb})$-morphism, which by induction
can only be a non-zero scalar multiple of the standard inclusion.
Hence $(x_{\alpha})^{m_{\alpha}}.\varphi(v_{w})$ is a non-zero multiple
of $v_\tau$. Further, by weight reasoning and ${\mathfrak{sl}}_2$-representation theory,
it follows that $x_{-\alpha}\varphi(v_w)=0$. By the usual exchange relation
we get
$$
x_{-\alpha}^{m_\alpha}(x_{\alpha})^{m_{\alpha}}\varphi(v_{w})
=c\varphi(v_w),
$$
for some nonzero complex number $c$, and hence $\varphi(v_w)$
is an extremal weight vector of weight $w(\Lam)$,
which finishes the proof.
\endpf
\begin{coro}\label{unique-morphism-coro}
Let $\tau < w$, then there exists (up to scalar multiples) a unique morphism of  $U(\tilde{\Lb})$-modules
$V_{\tau}(m \Lam_0) \longrightarrow V_{w}(m \Lam_0)$.
\end{coro}

Consider the Demazure module $D(m,n\Theta) = V_{-n\Theta}(m\Lam_0)$. We fix
a generator $w\neq0$ of the unique $U(\Lgc)$-fixed line in $D(m,\Theta)$.
Note (see \cite{FoL}) that $w$ spans the line of the highest weight vectors for $\Lhg$
in $V(m\Lam_0)$. By Theorem~\ref{demazure-zerlegung} we have
for $c_1\not=c_2$ an isomorphism
$$
D(m, (n+1) \Theta) \simeq D(m, \Theta)_{c_2} \ast D(m, n \Theta)_{c_1}.
$$
We extend this to an isomorphism of $U(\Lgc\oplus\bc K)$-modules
by letting $K$ operate on $D(m, (n+1) \Theta)$ by the level $m$,
and letting $K$ act on the second module by $0$ on the first factor
and on the second factor by the level $m$.  Define the map
$$
\tilde{\varphi}: D(m ,n\Theta)_{c_1} \longrightarrow D(m ,\Theta )_{c_2} \otimes D(m, n\Theta)_{c_1}
$$
by $\tilde{\varphi}(v) = w \otimes v$. This map is an $U(\Lgc)$-module morphism
because $w$ is $U(\Lgc)$-invariant, which extends, as above, to a
$U(\Lgc\oplus\bc K)$-module morphism.

The map respects the filtrations up to a shift:
let $v_2\in D(m ,\Theta )$ be a generator and
let $q$ be minimal such that $w\otimes v_2\in F^q(D(m ,\Theta )_{c_2} \otimes
D(m, n\Theta)_{c_1})$. By the $U(\Lgc)$-equivariance it follows that
$$
\tilde{\varphi}(F^j(D(m ,n\Theta )_{c_1})\subseteq
F^{j+q}(D(m ,\Theta )_{c_2} \otimes D(m, n\Theta)_{c_1})
$$
So we get an induced $U(\Lgc\oplus \bc K)$-morphism $\varphi$ between
the associated graded modules by $\varphi(\overline{v}) = \overline{ w\otimes v}$. 
$\varphi$ is nontrivial and so by Corollary \ref{unique-morphism-coro} it is
(up to multiplication by a scalar) the embedding of Demazure modules $\iota$. We proved:
\begin{lem} 
The map $\varphi: D(m, n\Theta) \longrightarrow D(m, \Theta)_{c_2} \ast D(m, n\Theta)_{c_1} \simeq 
D(m,(n+1)\Theta)$ induced by $\varphi(v) = \overline{w \otimes v}$ is an embedding of $U( \tilde{\Lb})$-modules.
\end{lem}

\noindent One knows that
$V(m\Lam_0) = \underset{n \rightarrow \infty}{\mbox{lim }}
D(m , n\Theta)$ as $U(\Lgc)$-modules, and also
as $U(\tilde{\Lb})$-modules. It follows by the above:
\begin{lem} Let $\Lg$ be a simple Lie algebra.\\
The following is a commutative diagram of $U(\Lgc)$-modules%
$$
\xymatrix{ D(m,n\Theta) \ar@{^{(}->}[rr]^\iota\ar[dd]_\wr& & D(m, (n+1)\Theta)\ar[dd]_\wr\\ \\
D(m,n\Theta)\ar[rr]^\varphi\ar[dd]_\wr & & D(m,\Theta)_{c_{n+1}} \ast D(m, n\Theta)_{c} \ar[dd]_\wr\\ \\
D(m, \Theta)_{c_n} \ast \ldots \ast D(m, \Theta)_{c_1}  \ar[rr]^\varphi&& D(m, \Theta)_{c_{n+1}} \ast 
D(m, \Theta)_{c_{n}} \ast \ldots \ast D(m,\Theta)_{c_1}\\
}
$$
where the down arrows are the isomorphism of Corollary~\ref{demazure-zerlegunghigh} %\ref{higherlevel-demazure}.
\end{lem}
\begin{thm}\label{limitthm}
Let $\Lg$ be a simple Lie algebra.
As $U(\Lgc)$-module, $V(m\Lambda_0)$ is isomorphic to the semi-infinite
fusion product
$$
V(m\Lambda_0)\simeq\underset{n \rightarrow \infty}{\mbox{lim }} 
D(m,\Theta) \ast \ldots \ast D(m,\Theta)
$$
\end{thm}
We expect the following to hold:
\begin{conj}\label{semiinfconjecture} 
Let $\Lam = m \Lam_0 + \lam$ be a dominant integral weight for $\Lhg$, then $V(\Lam)$ and
$$
\underset{n \rightarrow \infty}{\mbox{lim }} D(m,\Theta) \ast \ldots \ast D(m,\Theta) \ast V(\lam)
$$
are isomorphic as $\Lgc$-modules.
\end{conj}
\begin{rem}\rm
This isomorphism holds for the $\Lg$-module structure, see \cite{FoL}.
\end{rem}
\begin{rem}\label{moregenerallimit}
\rm 
As in \cite{FoL}, the limit construction above works in a much more general 
setting. Let $D(m ,\mu^{\vee})$ be a Demazure module with the property
that for some $k$ the fusion product 
$W=D(m ,\mu^{\vee})*\cdots*D(m ,\mu^{\vee}) \simeq D(m,k\mu^{\vee})$
% $W=D(m ,\mu^{\vee})^{\ast k}$
contains a highest weight vector of weight $m\Lam_0$.
Instead of $D(m,\Theta)$ one can then use the module $W$ in 
the direct limit construction above.
\end{rem}
%%%%%%%%%%%%%%%%%%%%%%%%%%%%%%%%%%%%%%%%%%%%%%%%%%%%%%%%%%%%%%%%%%
%%%%%%%%%%%%%%%%%%%%%%%%%%%%%%%%%%%%%%%%%%%%%%%%%%%%%%%%%%%%%%%%%%%%%
%%%%%%%% Bibliography %%%%%%%%%%%%%%%%%%%%%%%%%%%%%%%%%%%%%%%%%%%%%%%%%%%
%%%%%%%%%%%%%%%%%%%%%%%%%%%%%%%%%%%%%%%%%%%%%%%%%%%%%%%%%%
%%%%%%%%%%%%%%%%%%%%%%%%%%%%%%%%%%%%%%%%%%%%%%%%%%%%%%%%%%%%%%%%%%

\vskip 10pt \noindent
{{\eightsc  %\uppercase
{ G. Fourier, P. Littelmann: Mathematisches Institut der Universit\"at zu K\"oln,
Weyertal 86-90, 50931~K\"oln, Germany; email: gfourier@math.uni-koeln.de,
littelma@math.uni-koeln.de}}}
\end{document}